\documentclass[11pt]{amsart}
\usepackage{latexsym,amscd,xy,xspace}
\xyoption{v2}

\newcommand{\duplex}{\special{!statusdict begin true setduplexmode end}}

\newtheorem{Th}{Theorem}[section]
\newtheorem{Le}[Th]{Lemma}
\newtheorem{Co}[Th]{Corollary}
\newtheorem{Prop}[Th]{Proposition}

\theoremstyle{definition}

\newtheorem{Exa}{Example}[section]

\theoremstyle{remark}
\newtheorem{Rem}{Remark}[section]

\newcommand{\Proof}{\noindent{\bf Proof. \ }}

\renewcommand{\theenumi}{(\Alph{enumi})}

\newcommand{\n}{\noindent}
\newcommand{\HD}[2]{\ensuremath{\mathbb{H}^{#1}(M,\Z(#2)^{\infty}_D)}}
\newcommand{\HHD}[2]{\ensuremath{\mathbb{H}^{#1}(M,\Z(#2)_D)}}
\renewcommand{\H}[2]{\mbox{$H^{#1}(M,\mathbb{#2})$}}
\newcommand{\RBox}{$\hfill\Box$}
\newcommand{\coproduct}{\coprod\limits}
\newcommand{\cl}[1]{\ensuremath{\mathcal #1}}
\newcommand{\seq}[2][n]{\ensuremath{#2_0, \ldots , #2_{#1}}}

\renewcommand{\u}[1]{\ensuremath{\underline{#1}_M}}
\newcommand{\ux}[1]{\ensuremath{\underline{#1}_X}}
\newcommand{\us}[1]{\ensuremath{\underline{#1}^\ast_M}}
\newcommand{\cinf}[1]{\ensuremath{C^\infty(M, #1 )}}

\newcommand{\tC}{\ensuremath{\tilde{C}}}

\newcommand{\wt}{\widetilde}
\newcommand{\na}{\nabla}
\newcommand{\del}{\partial}
\newcommand{\db}{\bar{\partial}}
\newcommand{\bcs}[1]{\ensuremath{B^{#1}\mathbb{C}^\ast}}

\newcommand{\les}{long exact sequence\xspace}

\newcommand{\ses}{short exact sequence\xspace}
\newcommand{\sess}{short exact sequences\xspace}
\newcommand{\nbhd}{neighborhood\xspace}

\newcommand{\dlog}{\operatorname{dlog}}
\newcommand{\tors}{\operatorname{Tors}}
\newcommand{\tot}{\operatorname{Tot}}
\newcommand{\im}{\operatorname{im}}
\newcommand{\supp}{\operatorname{supp}}
\newcommand{\map}{\operatorname{Map}}

\newcommand{\C}{\ensuremath{\mathbb{C}}\xspace}
\newcommand{\N}{\ensuremath{\mathbb{N}}\xspace}
\newcommand{\Cst}{\ensuremath{\mathbb{C}^\ast}\xspace}
\newcommand{\D}{\ensuremath{\mathbb{D}}\xspace}
\renewcommand{\P}{\ensuremath{\mathbb{P}}\xspace}
\newcommand{\R}{\ensuremath{\mathbb{R}}\xspace}
\newcommand{\Z}{\ensuremath{\mathbb{Z}}\xspace}
\newcommand{\bH}{\ensuremath{\mathbb{H}}\xspace}

\newcommand{\cA}{\ensuremath{\mathcal{A}}\xspace}
\newcommand{\A}[1]{\ensuremath{{\mathcal A}^{#1}_{M,\mathbb{C}}}\xspace}

\newcommand{\ga}{\ensuremath{\alpha}\xspace}
\newcommand{\gb}{\ensuremath{\beta}\xspace}

\newcommand{\gG}{\ensuremath{\Gamma}\xspace}
\newcommand{\gd}{\ensuremath{\delta}\xspace}
\newcommand{\gD}{\ensuremath{\Delta}\xspace}

\newcommand{\gf}{\ensuremath{\varphi}\xspace}

\newcommand{\gL}{\ensuremath{\Lambda}\xspace}

\newcommand{\gs}{\ensuremath{\sigma}\xspace}

\newcommand{\gO}{\ensuremath{\Omega}\xspace}
\renewcommand{\go}{\ensuremath{\omega}\xspace}
\newcommand{\gp}{\ensuremath{\psi}\xspace}

\newcommand{\OO}[1]{\ensuremath{\Omega^{#1}_X}}
\newcommand{\Ost}{\ensuremath{\mathcal{O}^\ast_X}\xspace}

\newcommand{\gog}{\ensuremath{\mathbf{g}}\xspace}

\newcommand{\lra}{\longrightarrow}
\newcommand{\ra}{\rightarrow}
\newcommand{\hra}{\hookrightarrow}

\newcommand{\ras}{\mapsto}
\newcommand{\vl}{\bigl\lvert}
\newcommand{\vr}{\bigr\rvert}

\newcommand{\rh}{\rto|<\hole|<<\ahook}

\begin{document}
\duplex

\title[Geometry of Deligne cohomology]
{Geometry of  Deligne cohomology}
\author[Pawe\l \ Gajer]{Pawe\l \ Gajer}
\date{\today}
\address{Department of Mathematics, Texas A\&M University, College
Station, TX 77843-3368}
\email{gajer@math.tamu.edu}
\maketitle

The aim of this paper is to give a geometric interpretation of
holomorphic and smooth Deligne cohomology. Before stating the main
results we recall the definition and basic properties of Deligne
cohomology.

Let $X$ be a smooth complex projective variety and let $\OO{r}$ be the sheaf
of germs of holomorphic $r$-forms on $X$. The $q$th {\em
  Deligne complex of $X$} is the complex of sheaves
\[ \Z(q)_D :\quad \ux{\Z(q)} \lra \OO{0} \overset{d}{\lra} \OO{1}
\overset{d}{\lra} \cdots \overset{d}{\lra} \OO{q-1},
\]
where $\Z(q) = (2\pi \sqrt{-1})^q\Z \subset \C$, and $\ux{\Z(q)}$ is the
constant sheaf on $X$ associated with the group $\Z(q)$.  The
hypercohomology $\bH^{\ast}(X,\Z(q)_D)$ of the complex $\Z(q)_D$ is
called the {\em Deligne cohomology of $X$}. Our basic reference for
Deligne cohomology is \cite{esn&vie-dc}.

One of the key properties of Deligne cohomology is that for every
$p\geq 1$ the group $\bH^{2p}(X, \Z(p)_D)$ is the extension
\begin{gather}\label{basicholses}
0\lra J^p(X)\lra \bH^{2p}(X, \Z(p)_D) \lra H^{p,p}_{\Z}(X)\lra 0
\end{gather}
of the group $H^{p,p}_{\Z}(X)$ of integral $(p,p)$-classes of $X$ by
the the $p$th intermediate Jacobian $J^p(X)$ of Griffiths. For $p=1$
the group \HHD{2}{1} is isomorphic to the first cohomology group
$H^1(\Ost)$ of the sheaf \Ost of germs of non-vanishing holomorphic
functions on $X$, and the sequence \eqref{basicholses} reduces in this
case to the classical \ses
\begin{gather}\label{dim1ses}
0\lra J(X)\lra H^1(X, \Ost) \lra H^{1,1}_{\Z}(X)\lra 0
\end{gather}
It is well known that the group $H^1(X, \Ost)$ is isomorphic to the group
$CH^1(X)$
of divisors of $X$ modulo rational equivalence, the Jacobian
$J(X)$ is isomorphic to the group $CH_{hom}^1(X)$ of rational
equivalence classes of homologous to 0 divisors of $X$, and
$H^{1,1}_{\Z}(X)$ is the image of the cycle map $CH^1(X)\ra
H^2(X;\Z)$. One would like to have a similar cohomological description
of the groups $CH^p(X)$ of codimension $p$ cycles of $X$ modulo
rational equivalence, for $p>1$, together with an analogous to
\eqref{dim1ses} \ses completely describing the image
$H^{2p}_{alg}(X;\Z)$ and the kernel $CH_{hom}^p(X)$ of the cycle
homomorphism
\begin{gather}\label{cyclehom}
CH^p(X) \lra H^{p,p}_{\Z}(X).
\end{gather}
Deligne cohomology can be thought as a step in this direction. Indeed,
the cycle homomorphism \eqref{cyclehom} lifts to a homomorphism
\[CH^p(X) \lra H^{2p}(X;\Z(p)_D),\]
so that the diagram
\[{\diagram
  0\rto& CH_{hom}^p(X)\dto^{AJ} \rto& CH^p(X)\dto
  \rto& H^{2p}_{alg}(X;\Z)\dto \rto& 0\\ 0\rto&
  J^p(X)\rto& \bH^{2p}(X, \Z(p)_D) \rto& H^{p,p}_{\Z}(X)\rto&0
  \enddiagram}
\]
commutes, where $AJ: CH_{hom}^p(X) \ra J^p(X)$ is the Abel-Jacobi
homomorphism.

Another important property of Deligne cohomology is the existence
of characteristic classes, called ``regulators'',
\[c_{n,p} : K_n^{alg}(X) \lra \bH^{2p-n}(X;\Z(p)_D)\]
from the algebraic $K$-groups of $X$ into the Deligne cohomology of
$X$, which generalize the classical Chern classes of holomorphic
vector bundles.  Several important conjectures of arithmetic algebraic
geometry are formulated in terms of these regulators
\cite{bei-Lf,rss-beil}.

The second degree Deligne cohomology groups $\bH^2(X;\Z(q)_D)$ have
the following geometric interpretations.
\begin{itemize}
\item  $\bH^2(X;\Z(0)_D)$ is the ordinary second cohomology group
  $H^2(X;\Z)$ of $X$ that can be identified with the group of smooth
  principal $\Cst$-bundles over $X$.
\item  $\bH^2(X;\Z(1)_D)$ is isomorphic to the group $H^1(\Ost)$ of
  isomorphism classes of holomorphic principal $\Cst$-bundles over $X$.
\item $\bH^2(X;\Z(2)_D)$ is isomorphic to the group of isomorphism
  classes of holomorphic principal $\Cst$-bundles over $X$ with
  holomorphic connections.
\item For every $q>2$ the group $\bH^2(X;\Z(q)_D)$ is isomorphic to
  the group of isomorphism classes of holomorphic principal
  $\Cst$-bundles with flat connections over $X$.
\end{itemize}

Thanks to the above description of the groups $\bH^2(X;\Z(q)_D)$ the
geometric structure of regulators, cycle homomorphisms, and
Abel-Jacobi homomorphisms has been completely understood in the case
of divisors. It is expected that a geometric interpretation of higher
degree Deligne cohomology will lead to a better understanding of
regulators, cycle homomorphisms, and Abel-Jacobi homomorphisms for
cycles of codimension greater than one (see \cite{bry-holgerbs}).
J-L. Brylinski and P. Deligne found a geometric description of the
groups $\bH^3(X;\Z(3)_D)$ in terms of holomorphic gerbes with
connective structure and curving (see \cite{bry-greenbook}). The main
drawback of this description is that it very difficult to generalize
to higher degrees. In this paper we show that for every $q \geq 0$ and
$p\geq 2$ the Deligne cohomology groups $\bH^p(X;\Z(q)_D)$ and the
intermediate Jacobians $J^p(X)$ have geometric interpretations
naturally extending the geometric description of the groups
$\bH^2(X;\Z(q)_D)$ and the classical Jacobian $J(X)$ (see Theorem~D).
In particular, for every smooth complex projective variety $X$ we give
a geometric interpretation of the \ses \eqref{basicholses} that
generalizes the classical diagram.

{\small \[{\diagram 0\rto& J(X)\rto\dto^-\cong&
    H^1(X, \cl{O}^\ast_X) \rto\dto^-\cong&
    H^{1,1}_{\Z}(X)\rto\dto^-\cong&0\\ 0\rto& \left\{ \Text{iso.
        classes of\\ topologically\\ trivial\\ holomorphic\\
        $\C^{\ast}$-bundles\\ over $X$} \right\}\rto& \left\{
      \Text{iso. classes of\\ holomorphic\\ principal\\
        $\C^{\ast}$-bundles\\ over $X$} \right\}\rto& \left\{
      \Text{iso. classes of \\ smooth principal\\ $\C^{\ast}$-bundles
        over $X$\\ admitting\\ holomorphic\\ structures} \right\}
    \rto& 0 \enddiagram}\]}

The description of geometry of holomorphic Deligne cohomology will be
preceded by a discussion of a geometric interpretation of a smooth
counterpart of Deligne cohomology. The \emph{smooth Deligne
  cohomology} of a smooth manifold $M$ is the hypercohomology
$\HD{\ast}{q}$ of the $q$th {\em smooth Deligne complex of $M$}, which
is the complex of sheaves
\[ \Z(q)^{\infty}_D :\quad \u{\Z(q)} \lra \A{0} \overset{d}{\lra} \cdots
\overset{d}{\lra} \A{q-1},
\]
where $\cl{A}^n_M$ is the sheaf of germs of smooth differential
$n$-forms on $M$, and $\A{n} =\cl{A}^n_M\otimes \C$. Smooth Deligne
cohomology groups, similarly to odrinary Deligne cohomology groups,
are extensions of ordinary cohomology groups. Moreover, they can be
identified with Cheeger-Simons differential characters groups
\cite{che&sim-diffchar}. Our basic reference for smooth Deligne
cohomology and the pertinent homological algebra of sheaves is
\cite{bry-greenbook}.

Degree two smooth Deligne cohomology groups have similar to
holomorphic Deligne cohomology groups geometric interpretations with
smooth principal \Cst-bundles and smooth connections taking place of
holomorphic line bundles and holomorphic connections. The \ses
\eqref{basicholses} has the following two counterparts
\begin{gather}
  0 \lra \dfrac{A^{p-1}_{\C}(M)}{A^{p-1}_{\C}(M)_0} \lra \HD{p}{p}
  \lra \H{p}{Z} \lra 0 \label{1ses}\\
  \intertext{and for $p<q$}
  0 \lra \dfrac{\H{p-1}{C}}{\H{p-1}{Z}_{TF}} \lra \HD{p}{q} \lra
  \tors\H{p}{Z} \lra 0\label{2ses}
\end{gather}
where $A^{p-1}_{\C}(M)$ is the group of $\C$-valued (p-1)-forms on
$M$, $A^{p-1}_{\C}(M)_0$ is the subgroup of $A^{p-1}_{\C}(M)$
consisting of closed (p-1)-forms with integral periods, the group
$\H{p}{Z}_{TF}$ is the image of $\H{p}{Z}$ in $\H{p}{C}$, and
$\tors\H{p}{Z}$ is the torsion part of the group $\H{p}{Z}$. For $p=2$
these sequences have the following geometric interpretations.

\n{\bf Proposition~A. } {\em For every smooth manifold $M$ there is a
  commutative diagram} {\small
\[{\diagram
    {0} \rto & {\dfrac{A^{1}_{\C}(M)}{A^{1}_{\C}(M)_0}} \rto
    \dto^-{\cong}& {\HD{2}{2}} \rto \dto^-{\cong}& {\H{2}{Z}}
    \rto\dto^-{\cong}& {0}\\
    {0} \rto &
    {\left\{ \Text{iso. classes of\\connections on\\ $\C^{\ast}\times
          M \ra M$} \right\}}
    \rto & {\left\{ \Text{iso. classes of\\ smooth principal\\
          $\C^{\ast}$-bundles\\ with connections\\ over $M$} \right\}}
    \rto& {\left\{ \Text{iso. classes of\\ smooth principal\\
          $\C^{\ast}$-bundles\\ over $M$} \right\}} \rto& {0}
    \enddiagram}\]}
\n {\em with exact rows and the vertical arrows isomorphisms.}

\n{\bf Proposition~B.} {\em For every smooth manifold $M$ and every
  $q>2$ there exists a commutative diagram} {\small
\[{\diagram {0} \rto& {\dfrac{\H{1}{C}}{\H{1}{Z}_{TF}}} \rto
    \dto^-{\cong}& {\HD{2}{q}} \rto \dto^-{\cong}&
    {\text{Tors}\H{2}{Z}}\dto^-{\cong} \rto& {0} \\
{0} \rto& {\left\{\Text{iso. classes of\\flat connections\\ on
      $\C^{\ast}\times M \ra M$} \right\}} \rto & {\left\{ \Text{iso.
      classes of\\ smooth principal\\ $\C^{\ast}$-bundles with\\ flat
          connections\\ over $M$} \right\}} \rto&
    {\left\{\Text{iso.
          classes of\\ smooth principal\\ $\C^{\ast}$-bundles\\ over
          $M$\\ admitting\\ flat connections} \right\}} \rto& {0}
    \enddiagram}\]}
\n {\em with exact rows and the vertical arrows isomorphisms.}

Smooth Deligne cohomology is a natural framework for formulation of
the  Weil-Kostant Integrality Theorem (see Theorems~2.2.14
and 2.2.15 in \cite{bry-greenbook}), which is essentially equivalent
to the following result.

\n{\bf Proposition C. } {\em For every smooth manifold $M$ there is a
  commutative diagram} {\small
\[{\diagram
{0} \rto& {\H{1}{C^{\ast}}} \rto\dto^-{\cong}& {\HD{2}{2}}
\rto\dto^-{\cong}& {A^{2}_{\C}(M)_0} \dto^-{\cong} \rto& {0}\\
{0} \rto& {\left\{ \Text{iso. classes of\\ flat connections\\ on
      smooth\\ principal\\ $\C^{\ast}$-bundles\\ over $M$} \right\}}
\rto& {\left\{ \Text{iso. classes of\\ smooth principal\\
      $\C^{\ast}$-bundles with\\ connections\\ over $M$} \right\}}
\rto& {\left\{ \Text{curvatures of\\ connections \\on smooth\\
      principal \\ $\C^{\ast}$-bundles\\ over $M$} \right\}} \rto& {0}
\enddiagram}\]}
\n {\em with exact rows and the vertical arrows isomorphisms.}

Later on we will describe higher order analogues of the Weil-Kostant
Integrality Theorem.

Since $\Z(0)_D^\infty = \Z(0)_D = \u{\Z}$ both smooth and holomorphic
Deligne cohomology specialize to ordinary cohomology. Therefore, a
geometric interpretation of Deligne cohomology induces a geometric
interpretation of ordinary cohomology. Actually, it is natural to
start from a geometric description of ordinary cohomology, and then
enhance it to get a geometric model of smooth and holomorphic Deligne
cohomology. We proceed as follows.

The interpretation of \H{2}{Z} as the group of isomorphism classes of
smooth principal \Cst-bundles over $M$ comes from the isomorphism
\[      H^1(\us{\C}) \cong \H{2}{Z},    \]
given by the coboundary homomorphism in the cohomology \les associated
with the exponential \ses
\[\begin{CD} 0 @>>> \u{\Z} @>{\times
    2 \pi i}>> \u{\C} @>\exp>> \u{\C}^\ast @>>> 0 \end{CD}\] where for
any Lie group $G$ the symbol \u{G} stands for the sheaf of germs of
smooth G-valued functions on $M$.

A geometric interpretation of the groups \H{p}{Z}, for $p>2$, is
derived from a generalized exponential sequence, which is constructed
as follows. For every $s\geq 1$, the iterated classifying space of the
group \Cst
\[ B^{s}\Cst =
\underbrace{B(B(\cdots B}_{s \scriptstyle{\rm\ times}}(\Cst)\cdots))\]
can be equipped with a differentiable space structure so that the
homomorphisms in the \ses
\[0 \lra B^{s-1}\Cst \lra EB^{s-1}\Cst \lra B^{s}\Cst \lra 0\]
are smooth maps. The composition of these \sess induces an acyclic
resolution of the group \Z
\begin{gather*}
  0\lra \Z \lra \C \lra E\Cst \lra
  EB\Cst \lra EB^2\Cst \lra \cdots
  \lra EB^n\Cst \lra \cdots
\end{gather*}
which in turn induces an acyclic \emph{bar resolution} of the sheaf
\u{\Z}
\[0\lra \u{\Z}\lra \u{\C} \lra \us{E\C} \lra \us{EB\C} \lra
\us{EB^2\C} \lra \cdots \lra \us{EB^n\C} \lra \cdots
\]
where \us{EB^n\C} stands for the sheaf of {\em smooth}
$EB^n\Cst$-valued functions on $M$. The bar resolution of \u{\Z} is a
special case of a construction that assigns to a sheaf \cl{F} an
acyclic resolution of \cl{F}
\[0\lra \cl{F} \lra E\cl{F} \lra EB\cl{F}
\lra EB^2\cl{F} \lra \cdots
\lra EB^n\cl{F} \lra \cdots \]

The truncation of the bar resolution of the sheaf \u{\Z} in degree
$p-2$ induces the {\em generalized exponential \ses}\/
 {\small
\[\xymatrixnocompile{
  & &{\us{EB^{p-3}\C}}\rto& {\us{B^{p-2}\C}}\rto& {0} \\ & &\vdots
  \uto& & \\ & &{\us{EB\C}} \uto& & \\ & &{\us{E\C}} \uto& & \\
  {0}\rto &{\u{\Z}}\rto&{\u{\C}} \uto& & }
\]}
We show that the coboundary homomorphism
\[\gd : H^1(\us{B^{p-2}\C}) = H^{p-1}(\us{B^{p-2}\C}[-p+2])\lra
H^p(\u{\Z}) \] in the cohomology \les associated with the generalized
exponential sequence is an isomorphism. Since $H^1(\us{B^{p-2}\C})$ is
isomorphic to the group of isomorphism classes of {\em smooth}
principal $B^{p-2}\Cst$-bundles over $M$, we get the following result.

\textbf{Theorem H. } \emph{
  For every smooth manifold $M$ the group $H^{p}(M;\Z)$ is
  isomorphic to the group of isomorphism classes of smooth
  principal $B^{p-2}\Cst$-bundles over $M$.}


Theorem~H can be viewed as a generalization of J. Giraud's result
that identifies $H^3(M;\Z)$ with the group of equivalence classes of
gerbes bound by \us{\C} (see \cite{gir-nonAbCoh} and
\cite[Theorem~5.2.8]{bry-greenbook}). A correspondence between smooth
principal $B\Cst$-bundles and gerbes bound by \us{\C} is explained in
Appendix~A.

It is natural to expect that the groups \HD{p}{q} have a description
in terms of connections on smooth principal $B^{p-2}\Cst$-bundles over
$M$.  Indeed, we show that for every $s\geq 1$ the group $B^s\Cst$ is
equipped with a $B^s\C$-valued $B^s\Cst$-equivariant connection 1-form
that can be used to define connections on smooth principal
$B^s\Cst$-bundles.  It turns out that for $q>p$ the group $\HD{p}{q}$
is isomorphic to the group of isomorphism classes of smooth principal
$B^{p-2}\Cst$-bundles with flat connections over $M$.

A connection 1-form on a smooth principal $B^s\Cst$-bundle over $M$
has higher degree analogues called $k$-connections, $k=2, \ldots
,s+1$, which are defined inductively so that a relationship between
the $(k+1)$ and $k$-connections is analogous to the relation between a
connection on a principal $B^s\Cst$-bundle and a set of transition
functions of this bundle. The relation of isomorphism of principal
\Cst-bundles with connections can be generalized to an equivalence
relation on the set of $B^s\Cst$-bundles with $k$-connections such that
for every $p\geq 2$ the smooth Deligne cohomology group $\HD{p}{p}$ is
isomorphic to the group of equivalence classes of smooth principal
$B^{p-2}\Cst$-bundles with $k$-connections, $k=1,\dots,p-1$, over $M$.

The following two theorems give a geometric interpretation of the
\sess \eqref{1ses} and \eqref{2ses}.

\n{\bf Theorem A.} {\em For every smooth manifold $M$ and every $p\geq
  2$ there is a commutative diagram} {\small
  $${\diagram {0} \rto& {\dfrac{A^{p-1}_{\C}(M)}{A^{p-1}_{\C}(M)_0}}
    \rto\dto^-{\cong}& {\HD{p}{p}} \rto\dto^-{\cong}& {\H{p}{Z}}
    \rto\dto^-{\cong}& {0}\\ {0} \rto& {\left\{ \Text{equivalence
          classes\\ of $k$-connections\\ $k=1,2,\, \ldots \, ,p-1$\\
          on the trivial\\ B$^{p-2}\C^{\ast}$-bundle\\ over $M$}
      \right\}} \rto & {\left\{ \Text{equivalence classes\\ of smooth
          principal\\ B$^{p-2}\C^{\ast}$-bundles\\ with
          $k$-connections\\ $k=1,\dots,p-1$\\ over $M$} \right\}}
    \rto& {\left\{ \Text{iso. classes of \\ smooth principal\\
          B$^{p-2}\C^{\ast}$-bundles\\ over $M$} \right\}} \rto& {0}
    \enddiagram}$$} \n {\em with exact rows and the vertical arrows
  isomorphisms. }

Theorem~A generalizes the description due to J-L. Brylinski and P.
Deligne of the smooth Deligne cohomology group $\HD{3}{3}$ in terms of
equivalence classes of gerbes bound by $\u{\C}^\ast$ with connective
structures and curving (see \cite{bry-greenbook}). A procedure
associating with a connection on a smooth principal $B\Cst$-bundle
$E\ra M$ a connective structure on the associated with $E\ra M$ gerbe
is described in Appendix~A.


\n{\bf Theorem B.} {\em For every smooth manifold $M$ and every
  $q>p\geq 2$ there exists a commutative diagram}
{\small
  $${\diagram {0} \rto& {\dfrac{\H{p-1}{C}}{\H{p-1}{Z}_{TF}}}
    \rto\dto^-{\cong}& {\HD{p}{q}} \rto\dto^-{\cong}&
    {\text{Tors}\H{p}{Z}} \rto\dto^-{\cong}& {0}\\ {0} \rto& {\left\{
        \Text{iso. classes of\\ flat connections \\ on the trivial\\
          B$^{p-2}\C^{\ast}$-bundle\\ over $M$} \right\}} \rto &
    {\left\{ \Text{iso. classes of\\ smooth principal\\
          B$^{p-2}\C^{\ast}$-bundles\\ with\\ flat connections\\ over
          $M$} \right\}} \rto&
    {\left\{ \Text{iso. classes of\\ smooth
          principal\\ B$^{p-2}\C^{\ast}$-bundles\\ over $M$\\
          admitting\\ flat connections} \right\}} \rto& {0}
    \enddiagram}$$} \n {\em with exact rows and the vertical arrows
  isomorphisms. }

To an equivalence class $[E, \go_1,\, \ldots \, ,\go_{p-1}]$ of a
smooth principal $B^{p-2}\Cst$-bundle $E\ra M$ with $k$-connections
$-\go_1,\, \ldots \, ,(-1)^{p-1}\go_{p-1}$ one can assign a scalar
curvature
\[ s([E, \go_1,\, \ldots \, ,\go_{p-1}]) = (-1)^{p-1}d\go_{p-1},\]
which is a \C-valued differential $p$-form on $M$.  We show that, if
the scalar curvature of $[E, \go_1,\, \ldots \, ,\go_{p-1}]$ is zero,
then the sequence $(E, \go_1,\, \ldots \, ,\go_{p-1})$ is equivalent
to a unique up to isomorphism sequence $(E, \go', 0, \, \ldots \, ,
0)$ with $-\go'$ being a flat connection on $E\ra M$. This gives a
geometric interpretation of the \ses
\[ 0\lra \H{p-1}{C^{\ast}} \lra \HD{p}{p} \lra A^{p}_{\C}(M)_0 \lra 0.\]

\n{\bf Theorem C. } {\em For every smooth manifold $M$ and every
  $p\geq 2$ there is a commutative diagram} {\small
  \[{\diagram {0} \rto& {\H{p-1}{C^{\ast}}} \rto\dto^-{\cong}&
    {\HD{p}{p} } \rto\dto^-{\cong}& {A^{p}_{\C}(M)_0}
    \rto\dto^-{\cong}& {0}\\ {0} \rto& {\left\{ \Text{iso. classes
          of\\ flat connections\\ on\\ smooth principal\\
          B$^{p-2}\C^{\ast}$-bundles\\ over $M$} \right\}} \rto &
    {\left\{ \Text{equivalence classes\\ of smooth principal\\
          B$^{p-2}\C^{\ast}$-bundles\\ with $k$-connections\\
          $k=1,\dots,p-1$\\ over $M$} \right\}} \rto^-s& {\left\{
        \Text{scalar curvatures\\ of smooth principal\\
          B$^{p-2}\C^{\ast}$-bundles\\ over $M$} \right\}} \rto& {0}
    \enddiagram}\]} \n {\em with exact rows and the vertical arrows
  isomorphisms.}

The above diagram shows that scalar curvatures are closed forms with
integral periods, and that every closed form with integral periods is
a scalar curvature of a connection on a smooth principal
$B^s\Cst$-bundle. This generalizes the classical Weil-Kostant
Integrality Theorem.

If one considers in the place of the smooth Deligne complex
\[ \u{\Z(q)} \lra \u{\C} \overset{d}{\lra} \A{1}\overset{d}{\lra}\,
\cdots \, \overset{d}{\lra} \A{q-1}
\]
the complex
\[ \u{\Z(q)} \lra i\u{\R} \overset{d}{\lra} i\A{1}\overset{d}{\lra}\,
\cdots \, \overset{d}{\lra} i\A{q-1},
\]
where $i=\sqrt{-1}$, then the hypercohomology groups of the last
complex have essentially the same geometric interpretation as
\HD{p}{p}, with the only difference being that one has to replace
everywhere \Cst by the unit circle.

As was mentioned before the groups $\bH^2(X, \Z(q)_D)$ and $\bH^2(M,
\Z(q)_D^\infty)$ have similar geometric descriptions, with the only
difference being that $\bH^2(M,\Z(q)_D^\infty)$ is described in terms
of smooth principal \Cst-bundles and smooth connections and $\bH^2(X,
\Z(q)_D)$ is described in terms of holomorphic principal \Cst-bundles
and holomorphic connections. Exactly the same phenomenon takes place
in higher degrees. We define holomorphic principal $B^s\Cst$-bundles
and holomorphic $k$-connections on them and prove the following
holomorphic analogue of Theorem~A.

\pagebreak

\n \textbf{Theorem D. } \emph{For every smooth complex projective
  variety $X$ and every $p\geq 2$ and $q>0$ the group $\bH^r(X,
  \Z(q)_D)$ is isomorphic to the group of equivalence classes of
  holomorphic principal $B^{r-2}\Cst$-bundles over $X$ with
  holomorphic $k$-connections, for $k=1, 2, \ldots, q-1$. Moreover,
  there is a commutative diagram} {\small
  \[{\diagram
    0\rto& J^p(X)\rto\dto^-\cong& \bH^{2p}(X,\Z(p)_D )
    \rto\dto^-\cong& H^{p,p}_{\Z}(X)\rto\dto^-\cong &0 \\
   0\rto&
   \left\{ \Text{equivalence
          classes\\ of holomorphic\\
          $k$-connections\\ $k=1,2,\, \ldots \, ,p-1$\\
          on the topologically\\ trivial,  holomorphic\\
          B$^{2(p-1)}\C^{\ast}$-bundles\\ over $X$}
   \right\}\rto&
   \left\{ \Text{equivalence classes\\ of holomorphic\\
          principal\\ B$^{2(p-1)}\C^{\ast}$-bundles\\ with
          $k$-connections\\ $k=1,\dots,p-1$\\ over $X$}
   \right\}\rto&
   \left\{ \Text{iso. classes of \\ smooth principal\\
          B$^{2(p-1)}\C^{\ast}$-bundles\\ over $X$, admitting\\
          holomorphic\\ structures}
   \right\} \rto& 0  \enddiagram}\]}
\n\emph{ with exact rows and the vertical arrows isomorphisms. }

The paper is organized as follows. In Section~1 we define a
differentiable space structure on $B^s\Cst$ for every $s\geq 1$. In
Section~2 we discuss bar resolutions of sheaves. Section~3 is devoted
to smooth principal $B^sG$-bundles. There we show that for every
smooth manifold $M$ the group $H^p(M;\Z)$ is isomorphic to the group
of isomorphism classes of smooth principal $B^{p-2}\Cst$-bundles over
$M$. In Section~4 we study connections on smooth principal
$B^sG$-bundles and prove Theorem~B. In Section~5 we define
$k$-connections and prove Theorems A and C. In Section~6 we prove
Theorem~D. In Appendix~A we discuss correspondence between
\bcs{}-bundles with connections and gerbes with connective structures.
Appendix~B has been included for the convenience of the readers not
familiar with the geometric bar construction. In this appendix we
review basic properties of the geometric bar construction, and discuss
a geometric meaning of the relations appearing in the standard
definition of the classifying space $BG$.

\n {\bf Acknowledgments}\\ The inspiration to the my work on geometry
of Deligne cohomology came from lectures of Paulo Lima-Filho on
Deligne cohomology at Texas A\&M University.  Several conversations
with Paulo were very helpful in the early development of this work for
which I am grateful to him. An excellent introduction to the geometry
of Deligne cohomology is Jean-Luc Brylinski's book
\cite{bry-greenbook}. The examples and ideas contained in this book
served me as guiding principles in the studies on the structure of
smooth and holomorphic principal $B^s\Cst$-bundles.

\section{Differentiable structures on $EB^sG$ and $B^{s+1}G$}
In this section we define for any abelian Lie group $G$ a
differentiable space structure on the spaces $EB^sG$ and $B^{s+1}G$
for every $s\geq 1$. For the definition and basic properties of the
geometric bar construction we refer the reader to Appendix~B.

Let $M.$ be a simplicial smooth manifold. That is, $M.$ consists of a
family $\{M_n\}_{n\in \N}$ of smooth manifolds, together with smooth
maps
\[\del_i:M_n\lra M_{n-1},\quad s_i:M_n\lra M_{n+1},\]
where $i=0,1, \dots , n$, satisfying the identities
\begin{gather*}
  \del_i\del_j = \del_{j-1}\del_i \quad\text{for}\quad i<j,\\ s_is_j =
  s_{j+1}s_i \quad\text{for}\quad i\leq j,\\ \del_is_j =
\begin{cases}
  s_{j-1}\del_i \quad\text{for}\quad i<j,\\ \text{id}|_{M_n}
  \quad\text{for}\quad i=j,j+1\\ s_j\del_{i-1} \quad\text{for}\quad
  i>j+1
\end{cases}
\end{gather*}

The geometric realization $|M.|$ of $M.$ is the quotient space of the
disjoint union
\[\coprod_{n\geq 0} \gD^n\times M_n \]
with respect to the equivalence relation $\sim$ generated by the
relations
\begin{gather*}
  (\del^ix,m) \sim (x,\del_im) \quad \text{for}\quad (x,m)\in
  \gD^{n-1}\times M_{n},\\ (s^ix,m) \sim (x,s_im) \quad \text{for}\quad
  (x,m)\in \gD^{n+1}\times M_{n},
\end{gather*}
where the maps $\del^i : \gD^{n-1} \ra \gD^n$ and $s^i: \gD^{n+1} \ra
\gD^n$ are defined in the baricentric coordinates by
\begin{align*}
  \del^i(\seq[n-1]{x}) &= (\seq[i-1]{x},0,x_{i},\, \ldots \,
  ,x_{n-1}),\\ s^i(\seq[n+1]{x}) &=
  (\seq[i-1]{x},x_i+x_{i+1},x_{i+2},\, \ldots \, ,x_{n+1}).
\end{align*}

A {\em differentiable space structure} on the geometric realization
$|M.|$ of $M.$ consists of the class of all smooth \R-valued functions
on $|M.|$. We say that a function $f:|M.|\ra \R$ is {\em smooth} if
the composition
\[\begin{CD}
  \coproduct_{n \geq 0} \Delta^n \times M_n @>q>>|M.| @>f>> \R
\end{CD}\]
is smooth\footnote{A map $g:\Delta^n \times M_n\ra \R$ is smooth at
  $x\in \del\Delta^n \times M_n$ if there is an open neighborhood of
  $x$ in $H^n\times M_n$, where $H^n=\{(\seq{x})\in\R^{n+1}\, | \,
  \sum x_i=1\}$, and a smooth function $\tilde{g}:U\ra \R$ which
  restricted to $U\cap (\Delta^n \times M_n)$ coincides with $g$.},
where $q$ is the quotient map.  Equivalently, a smooth \R-valued
function on $|M.|$ is given by a family of smooth maps $f^n :\Delta^n
\times M_n \ra \R$ such that for every $n\geq 0$ the following two
diagrams commute
\[{\diagram
  {\gD^n\times M_n} \rto^-{f^n}& {\R} && {\gD^n\times M_n}
  \rto^-{f^n}& {\R} \\ {\gD^{n-1}\times M_n} \uto^{\del^i\times
    id}\rto_-{id\times\del_i}& {\gD^{n-1}\times M_{n-1}}
  \uto_{f^{n-1}}&& {\gD^{n+1}\times M_n} \uto^{s^i\times
    id}\rto_-{id\times s_i}& {\gD^{n+1}\times M_{n+1}} \uto_{f^{n+1}}&
  \enddiagram}\]

Let $M$ and $N$ are differentiable spaces (that is $M$ and $N$ are
spaces equipped with the appropriately defined classes of smooth
functions). Then a map $f:M\ra N$ is called a {\em smooth map} if for
every smooth function $g:N\ra \R$ the composition $g\circ f:M\ra \R$
is a smooth map.

It is easy to see that if $f. : M.\ra N.$ is a simplicial smooth map
between simplicial smooth manifolds $M.$ and $N.$, then $f.$ induces a
smooth map $|f.|: |M.|\ra |N.|$ between the geometric realizations of
$M.$ and $N.$ respectively.

\begin{Exa}\label{e1.2}
  With every Lie group $G$ there are associate simplicial smooth
  manifolds $G.$, $EG.$, and $BG.$ with simplicial smooth maps
\[ G. \lra EG. \lra BG. \]
whose geometric realizations give a universal principal $G$-bundle
\[ G \lra EG \lra BG \]
Thus, the inclusion $G \lra EG$ and the projection $EG \lra BG$ are
smooth maps.
\end{Exa}

{\Le \label{BV} Let $V$ be a vector space over a field $k$. Then $EV$
  and $BV$, taken with respect to the additive group structure of $V$,
  are $k$-vector spaces with respect to the following multiplication
  by scalars
  $$\begin{array}{ll} k\times EV \ra EV,& \hspace{.5cm} c \cdot \vl
    t_1, \cdots, t_n, v_0[v_1|\, \cdots \, |v_n]\vr = \vl t_1, \cdots,
    t_n, cv_0[cv_1|\, \cdots \, |cv_n]\vr \\ k\times BV \ra BV,&
    \hspace{.5cm} c \cdot \vl t_1, \cdots, t_n, [v_1|\, \cdots \,
    |v_n]\vr = \vl t_1, \cdots, t_n, [cv_1|\, \cdots \, |cv_n]\vr
\end{array}$$
Moreover, the projection $EV \ra BV$ is a linear map. }

The proof of Lemma~\ref{BV} is an easy exercise which we leave for the
reader.

\begin{Exa}\label{e1.0}
  Let $V$ be a separable \C-vector space. It is easy to see that the
  homomorphism
\[l:EV \lra V, \quad\quad l(|\seq{x},  \seq{v} |) =
\sum\limits_{i=0}^n x_iv_i\] is a splitting of the \ses
\[0\lra V \lra EV \lra BV \lra 0\]
We will show that $l:EV \lra V$ is a smooth map.

Let $\seq{e},\ldots$ be a base of $V$ and let $\pi_n:V\ra \C$ be the
projection on the subspace span by $e_n$. To prove smoothness of $l:EV
\lra V$ it is enough to show that for every $k\geq 0$ the composition
\[\begin{CD}
  EV @>l>> V @>\pi_k>> \C
\end{CD} \]
is smooth. But
\[ \pi_k(l(|\seq{x},  \seq{v} |)) =
\pi_k\bigl(\sum\limits_{i=0}^n x_iv_i\bigr) =
\bigl<\sum\limits_{i=0}^n x_iv_i, e_k \bigr> = \sum\limits_{i=0}^n
x_i<v_i,e_k>
\]
is a smooth map on $EV$. Hence $l:EV \lra V$ is smooth.
\end{Exa}

Let $G$ be an abelian Lie group. A differentiable space structure on
$EB^sG$ and $B^{s+1}G$, for $s\geq 1$, is defined by the following
inductive procedure.

Suppose we have a notion of a smooth function on $B^{s}G$ as well as
on each product $\gD^k\times (B^{s}G)^m$ for $k,m\geq 0$. Then
$f:EB^{s}G\ra \R$ is smooth if the composition
\[\begin{CD}
  \coproduct_{n \geq 0} \Delta^n \times (B^{s}G)^{n+1}@>q_E>> EB^{s}G
  @>f>> \R
\end{CD}\]
is smooth and $f:B^{s+1}G\ra \R$ is smooth if the composition
\[\begin{CD}
  \coproduct_{n \geq 0} \Delta^n \times (B^{s}G)^{n}@>q_B>> B^{s+1}G
  @>f>> \R
\end{CD}\]
is smooth. A function $f:\gD^k\times (B^{s+1}G)^m\ra \R$ is smooth if
the composition
\[\begin{CD}
  \gD^k\times(\coproduct_{n \geq 0} \Delta^n \times (B^{s}G)^{n})^m
  @>\text{id}\times (q_B)^m>> \gD^k\times (B^{s+1}G)^m @>f>> \R
\end{CD}\]
is smooth.

Directly from the above definition of differentiable structures on
$B^{s+1}G$ and $EB^sG$ it follows that all maps in the \ses
\[0\lra B^sG \lra EB^sG \lra B^{s+1}G\lra 0 \]
are smooth. It is also not difficult to see that the map
\[B^{s}G\times B^{s}G \lra B^{s}G,\quad (g,h)\ras gh^{-1} \]
is smooth.

A group $G$ carring a differentiable space structure so that the map
\[G\times G\ra G,\quad (g,h)\ras gh^{-1}\]
is smooth is called a {\em differentiable group}.

\begin{Exa}\label{e1.4}
 \label{p6} For every differentiable group $G$, there is a smooth deformational
 retraction $r: EG \times I \lra EG$ of $EG$ to $e \in EG$ which is a
 minor modification of the standard contraction from \cite{milg-bar}.
 In particular, if $G=B^s\Cst$, then for every $s\geq 1$ there is a
 smooth deformational retraction $r: EB^s\Cst \times I \lra EB^s\Cst$.

 The map $r$ is represented by the family of maps
\[r_n: (EG)_n \times I \lra (EG)_{n+1},
\]
where
\[ (EG)_n = q_E\bigl(\coproduct_{i\leq n} \Delta^i\times
G^{i+1}\bigr)\subset EG,
\]
and $r_n$ is defined by the formula
\[ r_n(|t_1,\dots, t_n, h_0[h_1|\dots |h_n]|, t) =
\vl \Phi(0,t), \Phi(t_1,t), \dots, \Phi(t_n,t), [h_0|h_1|\dots
|h_n]\vr,
\]
\n where $\Phi : [0,1]^2 \ra [0,1]$ is the composition
$$\Phi(x,t) = \phi(\min(1,x+t))$$ \n with $\phi : [0,1]\ra [0,1]$
being a smooth nondecreasing function so that $\phi(0) = 0$ and
$\phi(1) = 1$.

The contraction $r: EG \times I \lra EG$ is a smooth map, because for
every smooth function $g: EG \ra \R$ the diagram

$${\diagram {(\Delta^n\times G^{n+1})\times I}
  \rrto^-{\tilde{r}_n}\dto_{q_n \times id}& & {\Delta^{n+1}\times
    G^{n+2}} \dto^{q_E} \drto^{g^{n+1}} & \\ {EG \times I} \rrto_r& &
  {EG} \rto_g& {\R} \enddiagram}$$ \n commutes, where
$$\tilde{r}_n:(\Delta^n\times G^{n+1})\times I \lra \Delta^{n+1}\times
G^{n+2}$$ \n is a smooth map defined by the formula
$$\tilde{r}_n(t_1,\dots, t_n, h_0, h_1, \dots, h_n, t) =
\bigl(\Phi(0,t), \Phi(t_1,t), \dots, \Phi(t_n,t), e, h_0, h_1, \dots,
h_n\bigr).$$

\end{Exa}

\begin{Le}
  Suppose $M$ is a smooth manifold and $G$ is a differentiable group.
  If $f:M\ra BG$ is a map so that for every $x\in M$ there is an open
  neighborhood $U$ of $x$ in $M$ so that $f$ restricted to $U$ is of
  the form
\[f=\vl f_0, f_1, \ldots ,f_n, [g_1|\ldots |g_n]\vr     \]
with $f_0, f_1, \ldots ,f_n, g_1, \ldots ,g_n$ being smooth maps, then
$f$ is a smooth map.
\end{Le}

\Proof Suppose $g: BG\ra \R$ is a smooth map. Thus, for every $n\geq
1$ the composition
\[ \gD^n\times G^n \overset{q_B}{\lra}  BG \overset{g}{\lra} \R \]
is smooth. Consider the commutative diagram
\[{\diagram
  {M} \rrto^f\drto_{\bar{f}}&& {BG}\rto^g& {\R}\\ & {\gD^n\times
    G^n}\urto_{q_B}& \enddiagram}\] where $\bar{f}=(f_0, f_1, \ldots
,f_n, g_1, \ldots ,g_n)$. Since both $\bar{f}$ and $q_B\circ g$ are
smooth, the composition $f\circ g = \bar{f}\circ q_B\circ g$ is smooth
as well. Thus $f:M\ra BG$ is a smooth map. \qed

\section{Bar resolutions of sheaves}

The key to the geometric interpretations of the cohomology groups from
Theorems A, B, C, and D is the following construction of a bar resolution
of a sheaf.

Let $G$ be an abelian group. The composition of the \sess
\begin{gather}\label{basicses}
  0\lra B^nG \lra EB^nG \lra B^{n+1}G \lra 0
\end{gather}
\n induces the \les
\begin{gather}\label{basicles}
  0\lra G \lra EG \overset{\gs}{\lra} EBG \overset{\gs}{\lra} EB^2G
  \overset{\gs}{\lra} \, \cdots \, \overset{\gs}{\lra} EB^nG
  \overset{\gs}{\lra} \, \cdots \,
\end{gather}
\n where for every $n\geq 0$ the homomorphism
\[      \gs : EB^nG \lra EB^{n+1}G      \]
\n is the composition
\[        EB^nG \lra B^{n+1}G \lra EB^{n+1}G    \]
of the surjection $EB^nG \ra B^{n+1}G$ and the monomorphism $B^{n+1}G
\ra EB^{n+1}G$.

If $G$ is an abelian Lie (or differentiable) group, then, as we saw in
Example~\ref{e1.2}, the \ses \eqref{basicses} is a {\em smooth}
$B^sG$-extension of $B^{s+1}G$ (that is both $ B^sG \lra EB^sG$ and
$EB^sG \lra B^{s+1}G$ are smooth homomorphisms). Hence, the \les
\eqref{basicles} induces the \les of sheaves
\begin{gather}\label{basiclesofsh}
  0\lra \u{G} \lra \u{EG} \overset{\gs}{\lra} \u{EBG}
  \overset{\gs}{\lra} \u{EB^2G} \overset{\gs}{\lra} \, \cdots \,
  \overset{\gs}{\lra} \u{EB^nG} \overset{\gs}{\lra} \, \cdots \,
\end{gather}
which will be called the {\em bar resolution} of the sheaf $\u{G}$.

\begin{Prop}\label{acyclicprop}
  The sequence \eqref{basiclesofsh} is an acyclic resolution of the
  sheaf \u{G}.
\end{Prop}

\Proof It is enough to show that for every differentiable group $G$
the group $H^i(\u{EG})$ is trivial, for every $i>0$.

Recall, that a sheaf $\cl{F}$ on $X$ is {\em soft} if for every closed
subset $Z$ of $X$ the restriction map $\cl{F}(X) \ra \cl{F}(Z)$ is a
surjection. If $X$ is a paracompact space and $\cl{F}$ is a soft sheaf
on $X$, then $H^i(X;\cl{F})\cong 0$ for all $i>0$ (see
\cite[Theorem~1.4.6]{bry-greenbook} ).

{\Le\label{acyclic} For every differentiable group $G$ the sheaf
  $\u{EG}$ is soft.}

\Proof Let $Z$ be a close subset of $M$ and let $\sigma_Z$ be a
section of $\u{EG}$ over $Z$. By the definition of a section of a
sheaf over a closed set there is an open set $U \supset Z$ and an
extension $\sigma_U$ of $\sigma_Z$ to $U$. Since $M$ is paracompact,
there is a \nbhd $V$ of $Z$ such that $\bar{V} \subset U$. The
extension $\sigma$ of $\sigma_U$ (and hence also $\sigma_Z$) to a
global section of $\u{EG}$ is given by the formula
\[\sigma(x) = r(\sigma_U(x), \psi(x)),
\]
where $r: EG \times I \ra EG$ is the deformational retraction from
Example~\ref{e1.2}.B and $\psi: M \ra [0,1]$ is a smooth function
equal to 1 on $V$ and equal to 0 on $M - U$.  \qed

A bar resolution of the sheaf $\cA^k_M$ of germs of smooth differential
$k$-forms on $M$ is constructed as follows. Let $\gL^kT^\ast M$ be the
$k$th exterior power of the cotangent bundle $T^\ast M$ of $M$ and let
$E\gL^kT^\ast M$ and $B\gL^kT^\ast M$ be the associated with
$\gL^kT^\ast M$ bundles with fibers over $x\in M$ equal to
$E(\gL^kT^\ast_x M)$ and $B(\gL^kT^\ast_x M)$ respectively. The groups
$E(\gL^kT^\ast_x M)$ and $B(\gL^kT^\ast_x M)$ carry vector spaces
structures as in Lemma~\ref{BV}.

Let $E\cA^k_M$ and $B\cA^k_M$ be the sheaves of germs of smooth
sections of the vector bundles $E\gL^kT^\ast M$ and $B\gL^kT^\ast M$
respectively. A section \ga of the sheaf $E\cA^k_M$ over $U\subset M$
is of the form
\[ \ga = \vl f_0, \, \ldots \,  ,f_n, \ga_0, \, \ldots \, ,\ga_n\vr, \]
and a section $\gb$ of the sheaf $B\cA^k_M$ over $U$ is of the form
\[ \gb = \vl f_0, \, \ldots \,  ,f_n, [\gb_0:\; \cdots \; :\gb_n]\vr, \]
where $\ga_0, \, \ldots \, , \ga_n, \gb_0, \, \ldots \, , \gb_n$ are
smooth differential $k$-forms on $U$ and $\{f_i\}_{i=0}^n$ is a smooth
partition of unity on $U$. The group of section of the sheaf
$E\cA^k_M$ over an open set $U\subset M$ will be denoted by
$\gG(U,E\cA^k_M)$. Similarly, $\gG(U,B\cA^k_M)$ stands for the group
of sections of $B\cA^k_M$ over $U$.

Since the sequence of vector bundles
\[      0\lra \gL^kT^\ast M \lra E\gL^kT^\ast M \lra
B\gL^kT^\ast M \lra 0 \] is exact, the sequence of the groups
\[      0\lra \gG(U,\cA^k_M) \lra \gG(U,E\cA^k_M) \lra \gG(U,B\cA^k_M) \lra 0
\]
is exact, for every open subset $U$ of $M$. Hence, the sequence of
sheaves
\[      0\lra \cA^k_M \lra E\cA^k_M \lra B\cA^k_M \lra 0        \]
is exact. Similarly, if $EB^{s-1}\cA^k_M$ and $B^s\cA^k_M$ are the
sheaves of smooth sections of the vector bundles $EB^{s-1}\gL^kT^\ast
M$ and $B^s\gL^kT^\ast M$ respectively, then the sequence of sheaves
\[ 0\lra B^{s-1}\cA^k_M \lra EB^{s-1}\cA^k_M \lra B^s\cA^k_M \lra 0 \]
is exact. The composition of these sequences induces a \les
\begin{gather}\label{difffbarres}
  0\lra \cA^k_M \lra E\cA^k_M \overset{\gs}{\lra} EB\cA^k_M
  \overset{\gs}{\lra}\, \cdots \, \overset{\gs}{\lra} EB^s\cA^k_M
  \overset{\gs}{\lra} \, \cdots
\end{gather}
where
\[ \gs: EB^s\cA^k_M \lra EB^{s+1}\cA^k_M        \]
is the composition
\[      EB^s\cA^k_M \lra B^{s+1}\cA^k_M \lra EB^{s+1}\cA^k_M    \]
The sequence \eqref{difffbarres} will be called the {\em bar
  resolution} of the sheaf $\cA^k_M$.

A bar resolution of an arbitrary sheaf \cl{F} on a space $X$, which is
not necessarily a smooth manifold, can be defined as follows.

Let $E\cl{F}$ and $B\cl{F}$ be the sheaves associated with the
presheaves
\[ U \ras E(\cl{F}(U)) \qquad\text{and}\qquad U \ras B(\cl{F}(U))\]
respectively. Since the stalks of $E\cl{F}$ and $B\cl{F}$ at $x\in X$
are $E(\cl{F}_x)$ and $B(\cl{F}_x)$ respectively, where $\cl{F}_x$ is
the stalk of the sheaf \cl{F} at $x$, and the sequence
\[ 0\lra \cl{F}_x \lra E(\cl{F}_x) \lra B(\cl{F}_x) \lra 0      \]
is exact, the sequence of sheaves
\[ 0\lra \cl{F} \lra E\cl{F} \lra B\cl{F} \lra 0        \]
is exact. Iterating the above bar constructions we get for every
$s\geq 1$ the sheaves $EB^{s-1}\cl{F}$ and $B^s\cl{F}$ so that the
sequence
\[ 0\lra B^{s-1}\cl{F} \lra EB^{s-1}\cl{F} \lra B^s\cl{F} \lra 0        \]
is exact. The composition of these sequences gives the {\em bar
  resolution of \cl{F}}
\begin{gather*}
  0\lra \cl{F} \lra E\cl{F} \overset{\gs}{\lra} EB\cl{F}
  \overset{\gs}{\lra} EB^2\cl{F} \overset{\gs}{\lra} \, \cdots \,
  \overset{\gs}{\lra} EB^n\cl{F} \overset{\gs}{\lra} \, \cdots \,
\end{gather*}
The complex of sheaves
\begin{gather*}
  \cl{B}^\ast(\cl{F}):\qquad E\cl{F} \overset{\gs}{\lra} EB\cl{F}
  \overset{\gs}{\lra} EB^2\cl{F} \overset{\gs}{\lra} \, \cdots \,
  \overset{\gs}{\lra} EB^n\cl{F} \overset{\gs}{\lra} \, \cdots \,
\end{gather*}
will be called the {\em bar complex of \cl{F}}.

An easy modification of the proof of Lemma~\ref{acyclic} shows that
the bar resolution of \cl{F} is an acyclic resolution of \cl{F}.
Therefore, the cohomology of \cl{F} is equal to the cohomology of the
cochain complex
\begin{gather*}
  \gG(M,E\cl{F}) \overset{\gs}{\lra} \gG(M,EB\cl{F})
  \overset{\gs}{\lra} \gG(M,EB^2\cl{F}) \overset{\gs}{\lra} \, \cdots
  \, \overset{\gs}{\lra} \gG(M,EB^n\cl{F}) \overset{\gs}{\lra} \,
  \cdots \,
\end{gather*}
The above complex will be called the {\em bar cochain complex of
  \cl{F}} and we will denote it by $C_B^\ast(\cl{F})$.

Note, that the above construction applied to \u{G} and $\cA^{k}_M$
produces resolutions of \u{G} and $\cA^{k}_M$ that do not coincide
with the resolutions \eqref{basiclesofsh} and \eqref{difffbarres}.  In
the sequel, when referring to bar resolutions of \u{G} or $\cA^{k}_M$
we will always mean the resolutions \eqref{basiclesofsh} or
\eqref{difffbarres} respectively.

 The bar cochain complex $C_B^\ast(\u{G})$ of the sheaf \u{G} is
  of the form
\[\cinf{EG}\overset{\gs}{\lra}
\cinf{EBG}\overset{\gs}{\lra}\, \cdots \, \overset{\gs}{\lra}
\cinf{EB^nG}\overset{\gs}{\lra}\, \dots \] \n where
\[ \gs : \cinf{EB^nG}\lra  \cinf{EB^{n+1}G}\]
is the composition
\[ \cinf{EB^nG}\overset{\pi_\ast}{\lra}
\cinf{B^{n+1}G}\overset{i_\ast}{\lra} \cinf{EB^{n+1}G} \] \n with
$\pi: E(B^nG) \ra B(B^nG)=B^{n+1}G$ being the projection map of the
universal principal $B^nG$-bundle and $i: B^{n+1}G \ra EB^{n+1}G$
being the inclusion of the fiber into the total space of the universal
principal $B^{n+1}G$-bundle.

In a sense, the bar cochain complex of the sheaf \u{\Z}
\[\cinf{E\Z}\overset{\gs}{\lra}
\cinf{EB\Z}\overset{\gs}{\lra}\, \cdots \, \overset{\gs}{\lra}
\cinf{EB^n\Z}\overset{\gs}{\lra}\, \dots \] can be thought of as a
smooth version of Karoubi's complex (see \cite{kar-ftnc})
\[\map(X, AG(\D^1))\overset{\gs_\ast}{\lra} \map(X, AG(\D^2))
\overset{\gs_\ast}{\lra} \, \cdots \, \overset{\gs_\ast}{\lra} \map(X,
AG(\D^n)) \overset{\gs_\ast}{\lra} \, \cdots \] of the topological
non-commutative differential forms on a space $X$, where $AG(\D^n)$ is
the free abelian group on the disk $\D^n$, and
\[\gs_\ast : \map(X, AG(\D^n)) \lra \map(X, AG(\D^{n+1}))\]
is a homomorphism induced by the composition of maps
\[ \D^n \lra \D^n/\del\D^n = S^n = \del \D^{n+1} \hra \D^{n+1}. \]

{}From the functoriality of the geometric bar construction it follows
that the bar resolution of sheaves is functorial as well. Moreover,
since for every \ses of topological groups
\[ 0\lra K \lra G \lra H \lra 0 \]
the sequences
\begin{gather*}
  0\lra EK \lra EG \lra EH \lra 0\\ \intertext{and} 0\lra BK \lra BG
  \lra BH \lra 0
\end{gather*}
are exact, every \ses of sheaves
\[ 0\lra \cl{E} \lra \cl{F} \lra \cl{G} \lra 0  \]
induces a \ses of complexes of sheaves
\[ 0\lra \cl{B}^\ast(\cl{E}) \lra \cl{B}^\ast(\cl{F})
\lra \cl{B}^\ast(\cl{G}) \lra 0 \] Hence, every complex of sheaves
$\cl{F}^\ast$ has an acyclic resolution given by the total complex
$\tot^\ast(\cl{B}^\ast(\cl{F}^\ast))$ associated with the double
complex $\cl{B}^\ast(\cl{F}^\ast)$. The cohomology of $\cl{F}^\ast$ is
equal to the cohomology of the total cochain complex
$\tot^\ast(C_B^\ast(\cl{F}^\ast))$.

\begin{Exa}
  The double cochain complex $C_B^\ast(\cA^\ast_M)$ of the de Rham
  complex
\[\cA^\ast_M:\qquad \cA^0_M \overset{d}{\lra} \cA^1_M
\overset{d}{\lra} \cA^2_M \overset{d}{\lra} \, \cdots \,
\overset{d}{\lra} \cA^n_M \overset{d}{\lra} \cdots\] is given by the
diagram {\small
\[{\diagram
  & & &\\ {\gG(M, E\cA^2_M)}\uto^-d\rto^\gs& & & \\ {\gG(M,
    E\cA^1_M)}\uto^-d\rto^\gs& {\gG(M,EB\cA^1_M)}\uto^-d\rto^-\gs& \\
  {\gG(M,E\cA^0_M)}\uto^-d\rto^\gs& {\gG(M,EB\cA^0_M)}\uto^-d\rto^\gs&
  {\gG(M,EB^2\cA^0_M)}\uto^-d\rto^-\gs &\\ \enddiagram} \]} Note that
the $n$th column of this double complex is the complex of global
sections of the acyclic resolution
\[EB^{n-1}\cA^0_M\ra EB^{n-1}\cA^1_M\ra \cdots \ra EB^{n-1}\cA^s_M\ra
\cdots\] of the sheaf $\u{EB^{n-1}\R}^\gd$, where $\R^\gd$ is the
group \R taken with the discrete topology. Therefore, the total
complex of $C_B^\ast(\cA^\ast_M)$ is an acyclic resolution of the de
Rham complex of $M$ and the bar complex
\[\u{E\R}^\gd \lra \u{EB\R}^\gd \lra \u{EB^2\R}^\gd \lra\cdots \lra
\u{EB^s\R}^\gd\lra \cdots
\]
of $\u{\R}^\gd$. Thus, the bar complex of the de Rham complex
$\cl{A}^\ast_M$ of $M$ plays a similar role to the \v{C}ech complex of
$\cl{A}^\ast_M$ inducing an isomorphis between the de Rham cohomology
of $M$ and the sheaf cohomology of the constant sheaf $\u{\R}^\gd$.
\end{Exa}

\begin{Rem}
  There is a close relationship between the bar and \v{C}ech cochain
  complexes of the sheaf $\cA^k_M$. Actually, every smooth partition
  of unity $\{f_i\}_{i\in I}$ subordinated to an open covering
  $\cl{U}=\{U_i\}_{i\in I}$ of a manifold $M$ induces a cochain
  homomorphism
\[ \gf^\ast:\check{C}^\ast(\cl{U},\cA^k_M) \lra C_B^\ast(\cA^k_M)       \]
so that
\[ \gf^p : \check{C}^p(\cl{U},\cA^k_M) \lra C_B^p(\cA^k_M)=
EB^p\cA^k_M(M)\] is the composition
\[ \check{C}^p(\cl{U},\cA^k_M)\lra \check{C}^{p-1}(\cl{U},B\cA^k_M)\lra
\, \cdots \, \lra \check{C}^0(\cl{U},B^p\cA^k_M)\lra EB^p\cA^k_M(M)\]
where for $r>0$ the homomorphism
\[ \gf^{r,s} : \check{C}^r (\cl{U}, B^{s}\cA^k_M)
\lra \check{C}^{r-1} (\cl{U}, B^{s+1}\cA^k_M) \] is defined for $\xi =
\{ \xi_{\seq[r]{i}} \in B^s\cA^k_M(\bigcap\limits_{j=0}^r U_{i_j}) \}$
by the formula
\[      \gf^{r,s}(\xi)_{\seq[r-1]{i}} =
\vl f_{l_0}, \, \ldots \, ,f_{l_n}, [\xi_{l_0,\seq[r-1]{i}}: \, \cdots
\, :\xi_{l_n,\seq[r-1]{i}}]\vr \] and the homomorphism
\[ \gf^{0,p}:\check{C}^0(\cl{U},B^p\cA^k_M)\lra EB^p\cA^k_M(M)\]
is given by
\[      \gf^{0,p}(\{\xi_i\}) =
\vl f_{l_0}, \, \ldots \, ,f_{l_n}, [\xi_{l_0}:\, \cdots \, :
\xi_{l_n}]\vr. \]

\end{Rem}

\section{Smooth principal $B^s\C^{\ast}$-bundles}
In this section we will show that if $M$ is a smooth manifold, then
the group $H^k(M;\Z)$ can be identified with the group of isomorphism
classes of smooth principal $B^{k-2}S^1, B^{k-2}\Cst$, or $B^{k-1}\Z$
bundles over $M$.

Let $G$ be an abelian Lie group. A principal $B^sG$-bundle $E \ra M$
over a smooth manifold $M$ is {\em smooth} if the transition functions
of this bundle are smooth.

The proof of the following proposition, essentially due to tom Dieck
\cite{tDie-klass}, shows an explicit formula for a classifying map of
a smooth principal bundle in terms of its transition functions.

{\Prop \label{smoothpr} Let $G$ be a differentiable group. Then for
  every smooth principal $G$-bundle $\pi :E \ra M$ there is a smooth
  map $\gf:M \ra BG$ such that $E \ra M$ is the pull-back of the
  universal principal $G$-bundle by $\gf$. }

{\Co \label{smoothco} For every smooth principal $B^s\Cst$-bundle $E
  \ra M$ there is a smooth map $\gf:M \ra B^{s+1}\Cst$ such that $E
  \ra M$ is the pull-back of the universal principal $B^s\Cst$-bundle
  by $\gf$. }

\n{\bf Proof of Proposition~\ref{smoothpr}.} Let $\cl{U}
=\{U_i\}_{i\in I}$ be an open covering of $M$ so that for every $i\in
I$ there is a trivialization $$\gp_i : \pi^{-1}(U_i) \lra U_i\times
G.$$

Define $g_i : E \ra G$ by the formula
$$g_i(x)=
\begin{cases}
  \text{pr}_2(\gp_i(x))& \text{for}\quad x\in \pi^{-1}(U_i)\\ e &
  \text{for}\quad x\notin \pi^{-1}(U_i)
\end{cases}$$
\n where $e$ is the neutral element of $G$ and $\text{pr}_2 :
U_i\times G \ra G$ is the projection on the second factor.

Let $\{f_i\}_{i\in I}$ be a partition of unity subordinated to the
covering $\cl{U}$ and let \linebreak $\bar{\gf}:E \ra EG$ be the map
$$\bar{\gf}(y) = \vl f_{i_0}(\pi(y)), f_{i_1}(\pi(y)), \, \ldots \, ,
f_{i_n}(\pi(y)), g_{i_0}(y), g_{i_1}(y), \, \ldots \, ,
g_{i_n}(y)\vr,$$ \n where $\seq{i}$ are the indices so that for each
$i\in \{\seq{i}\}\quad f_i(\pi(y))\neq 0$. It is easy to see that
$\bar{\gf}$ is a $G$-equivariant map and hence it induces a morphism
of principal $G$-bundles
$${\diagram {E} \rto^{\bar{\gf}}\dto_{\pi}& {EG}\dto \\ {M}
  \rto_{\gf}& {BG} \enddiagram}$$ \n where the restriction of $\gf : M
\ra BG$ to $U_j \subset M$ is given by the formula
$$\gf (x) = \vl f_{i_0}(x), f_{i_1}(x), \, \ldots \, , f_{i_n}(x),
[g_{i_0}(\gs(x)): g_{i_1}(\gs(x)): \, \cdots \, :
g_{i_n}(\gs(x))]\vr,$$ \n where $\gs: U_j\ra \pi^{-1}(U_j)$ is a
smooth section of the restriction $\pi^{-1}(U_j) \ra U_j$ of $\pi:E\ra
M$ to $\pi^{-1}(U_j)$. Note that $\gf(x)$ does not depend on the
choice of the section $\gs$ because $g_i$s are $G$-equivariant maps.
In the non-homogeneous coordinates
\[      \gf (x) = \vl f_{i_0}(x), f_{i_1}(x),\, \ldots \,  ,f_{i_{n}}(x),
[g_{i_0i_1}(x)|g_{i_1i_2}(x)|\, \cdots \,
|g_{i_{n-1}i_n}(x)]\vr,
\]
where
\[      g_{ij}(x) = (g_i(\gs(x)))^{-1}\cdot g_j(\gs(x))        \]
are the transition functions of the bundle $E\ra M$ associated with
the open covering of $M$ by the sets $\{x\in M\, | \, f_i(x)>0\}$.
Since $g_i$ is smooth on $\supp (f_i)\cap U_j$ for every $i,j\in I$
and \gs is smooth on $U_j$, the map \gf is smooth on $U_j$ for every
$j\in I$ and hence \gf is smooth on $M$. T. tom Dieck showed in
\cite{tDie-klass} that $\gf : M \ra BG$ is the classifying map of the
bundle $\pi : E \ra M$ (tom Dieck works in the setting of Milnor's bar
construction, but all he does extends easily to the context of
Milgram's bar construction).  \qed

\begin{Exa}\label{e1.1}
  The isomorphism $H^2(\R^3 - 0,\Z) \cong H^1(\R^3 -
  0,\underline{\Cst}_{\R^3 - 0})$ implies that every element of
  $H^2(\R^3 - 0,\Z)$ corresponds to a unique isomorphism class of a
  smooth principal \Cst-bundle over $\R^3 - 0$. Let $L$ be a smooth
  principal \Cst-bundle over $\R^3 - 0$ representing a generator of
  $H^2(\R^3 - 0,\Z) \cong \Z$. The proof of Proposition~\ref{smoothpr}
  shows how to describe a smooth classifying map $\gp_L:(\R^3 -0) \ra
  B\Cst$ of $L$ in terms of some transition functions of $L$. Let $S^3
  = \R^3 \cup \{\infty\}$ and consider the open subsets $U_0 =\R^3$
  and $U_\infty = S^3-\{0\}$ of $S^3$. Since $U_0\cap U_\infty = \R^3
  - 0$ we can think of the classifying map $\gp_L:(\R^3 -0) \ra B\Cst$
  as a transition function of a smooth principal $B\Cst$-bundle $BL$
  over $S^3$. From the proof of Proposition~\ref{H(M,Z)} it follows
  that the isomorphism class of $BL$ corresponds to the generator of
  $H^3(S^3,\Z)$. Let $\widetilde{BL}$ be the pull-back of $BL$ by the
  standard retraction $(\R^4 - 0) \ra S^3$. The classifying map of the
  bundle $\widetilde{BL}$ can be identified with a transition function
  of a smooth principal $B^2\Cst$-bundle over $S^4$, representing a
  generator of $H^4(S^4,\Z)$. Iterating the above procedure we get a
  family of smooth principal $B^k\Cst$-bundle over $S^k$, representing
  generators of the groups $H^k(S^k,\Z)$ for $k\geq 2$.
\end{Exa}

\begin{Prop} \label{H(M,Z)} Let $G$ be one of the group $S^1, \Cst$,
  or $B\Z$. Then for every smooth manifold $M$ and every $p\geq 2$ the
  group $H^p(M,\Z)$\ is isomorphic to:
  \renewcommand{\theenumi}{(\roman{enumi})}
\begin{enumerate}
\item the group $L(B^{p-2}G)_M$ of isomorphism classes of {\em smooth}
  principal $B^{p-2}G$-bundles over $M$.
\item the group $[M, B^{p-1}G]^\infty$ of smooth homotopy classes of
  smooth maps from $M$ to $B^{p-1}G$.
\end{enumerate}
\end{Prop}

\n{\bf Proof of part (i) of Proposition~\ref{H(M,Z)}.} Since for any
abelian differentiable group $G$ the group of isomorphism classes of
smooth principal $G$-bundles over $M$ is isomorphic to $H^1(\u{G})$,
we have to prove that there is an isomorphism
\[      H^p(M;\Z) \cong H^1(\u{B^{p-2}G}).              \]

Consider the cohomology \les
\begin{gather*}
  \ra \bH^{p-1}(\u{EB^{<p-2}G})\ra H^1(\u{B^{p-2}G})\ra H^p(M;\Z)\ra
  \bH^{p}(\u{EB^{<p-2}G})\ra
\end{gather*}
associated with the generalized exponential sequence
\begin{gather}\label{expseq1}
  0\lra \u{\Z} \lra \u{EB^{<p-2}G} \lra \u{B^{p-2}G}[-p+2] \lra 0
\end{gather}
where \u{EB^{<p-2}G} is the complex
\begin{gather*}
  \u{H} \overset{\gs}{\lra} \u{EG} \overset{\gs}{\lra} \u{EBG}
  \overset{\gs}{\lra} \u{EB^2G} \overset{\gs}{\lra} \, \cdots \,
  \overset{\gs}{\lra} \u{EB^{p-3}G}
\end{gather*}
with \u{H} being equal to $\R, \C$, or $E\Z$ for $G=S^1, \Cst$, or
$B\Z$ respectively.

Since for every $s\geq 0$ the sheaf \u{EB^sG} is acyclic, the
cohomology of the complex \u{EB^{<p-2}G} is equal to the cohomology of
the cochain complex
\begin{gather*}
  H \overset{\gs}{\lra} \u{EG}(M) \overset{\gs}{\lra} \u{EBG}(M)
  \overset{\gs}{\lra} \u{EB^2G}(M) \overset{\gs}{\lra} \, \cdots \,
  \overset{\gs}{\lra} \u{EB^{p-3}G}(M)
\end{gather*}
of the groups of global sections of the components of \u{EB^{<p-2}G}.
Therefore, for every $q>p-2$
\[ \bH^{q}(\u{EB^{<p-2}G})\cong H^q(\u{EB^{<p-2}G}(M))\cong 0.\]
Hence, the coboundary homomorphism
\[H^1(\u{B^{p-2}G}) \lra H^p(M;\Z)      \]
in the cohomology \les associated with \eqref{expseq1} is an
isomorphism. \qed

\begin{Rem}\label{rem1}
  Let $G$ be an arbitrary abelian Lie group. Replacing in the proof of
  Proposition~\ref{H(M,Z)} the sequence \eqref{expseq1} be the
  appropriate \ses associated with the bar resolution of \u{G}, we
  would get an isomorphism between $H^p(\u{G})$ and the group
  $L(B^{p-2}G)_M$ of isomorphism classes of {\em smooth} principal
  $B^{p-2}G$-bundles over $M$.
\end{Rem}

Part (ii) of Proposition~\ref{H(M,Z)} is a straightforward consequence
of the following lemma.

\begin{Le}\label{classification1}
  Let $G$ be a differentiable group. Then the group of isomorphism
  classes of smooth principal $G$-bundles over $M$ is isomorphic to
  the group $[M, BG]^\infty$ of smooth homotopy classes of smooth maps
  from $M$ to $BG$.
\end{Le}

\n{\bf Proof. } Let $G$ be a differentiable group. We will show that
there is an isomorphism
\[      [M, BG]^\infty \cong H^1(\u{G}).                \]
The beginning of the cohomology \les associated with the \ses
\[ 0 \lra \u{G} \lra \u{EG} \lra \u{BG} \lra 0  \]
\n is of the form
\[      \, \cdots \,  \lra C^\infty(M,EG) \overset{\pi_\ast}{\lra}
C^\infty(M,BG) \lra H^1(\u{G}) \lra H^1(\u{EG}) \lra \, \cdots  \]

Since $H^1(\u{EG}) \cong 0$, we have the isomorphism
\[      \dfrac{C^\infty(M,BG)}{\pi_\ast C^\infty(M,EG)} \cong
H^1(\u{G}). \] The image $\pi_\ast C^\infty(M,EG)$ of the group
$C^\infty(M,EG)$ in $C^\infty(M,BG)$ consists of those smooth maps
from $M$ to $BG$ that lift to maps from $M$ to $EG$. It is easy to see
that $f:M\ra BG$ has a lift to $\tilde{f}: M \ra EG$ if and only if
$f$ is smooth homotopic to a constant map. Hence
\[      \dfrac{C^\infty(M,BG)}{\pi_\ast C^\infty(M,EG)} \cong [M,
BG]^\infty      \] \qed

\begin{Rem}
  Let $G$ be a topological group. Replacing in the proof of
  Lemma~\ref{classification1} the sheaves of smooth maps on $M$ by
  sheaves of continues maps on some space $X$, we get an isomorphism
  between the group of isomorphism classes of principal $G$-bundles
  over $X$ and the group $[X,BG]$ of homotopy classes of maps from $X$
  to $BG$.
\end{Rem}

\section{Flat Connections on principal $B^s\C^{\ast}$-bundles}

In this section we show that for every $s\geq 1$ the group
$B^s\C^{\ast}$ is equipped with the canonical
$B^s\C^{\ast}$-equivariant $B^s\C$-valued connection 1-form
$B^s(z^{-1}dz)$. By analogy with the Lie group case, the form
$B^s(z^{-1}dz)$ is used to define connections on smooth principal
$B^s\C^{\ast}$-bundles. We show that for $q>p$ the smooth Deligne
cohomology group $\HD{p}{q}$ is isomorphic to the group of isomorphism
classes of smooth principal $B^{p-2}\C^{\ast}$-bundles with flat
connections. Moreover, we prove Theorem~B.

\subsection{The canonical connection 1-forms on $B^s\C^{\ast}$ and
  $EB^s\C^{\ast}$}

Let $M.$ be a simplicial smooth manifold. A {\em smooth $p$-form}
$\alpha$ on the geometric realization $|M.|$ of $M.$ is a family
$\{\alpha^n\}$ of differential $p$-forms $\alpha^n$ on $\Delta^n\times
M_n$ satisfying for every $0\leq i \leq n$ the following compatibility
conditions
\begin{gather}
  (\del^i\times id )^{\ast}\alpha^n = (id \times \del_i
  )^{\ast}\alpha^{n-1}\label{compcond1}\\ (s^i\times id
  )^{\ast}\alpha^n = (id \times s_i
  )^{\ast}\alpha^{n+1}\label{compcond2}
\end{gather}
\n where $\del^i\times id,\, id \times \del_i,\, s^i\times id$, and
$id \times s_i$ are the maps
$${\diagram {\gD^{n-1}\times M_{n-1}} & {\gD^{n-1}\times M_n}
  \lto_-{id\times \del_i}\rto^{\del^i\times id}& {\gD^n\times M_n}\\
  {\gD^{n+1}\times M_{n+1}} & {\gD^{n+1}\times M_n} \lto_-{id\times
    s_i}\rto^{s^i\times id}& {\gD^n\times M_n} \enddiagram}$$ \n with
$\del^i$ and $s^i$ being the coface and the codegeneracy maps on
$\Delta^n$s and $\del_i$, $s_i$ being the face and the degeneracy maps
on $M_n$s.

\begin{Exa} \label{e2.1}\end{Exa}
\begin{enumerate}
\item Let $G$ be a Lie group and let $g^{-1}dg$ be the canonical
  \gog-valued connection 1-form on $G$, where \gog is the Lie algebra
  of $G$. The total space $EG$ of the universal principal $G$-bundle
  $EG\ra BG$ carries a smooth \gog-valued form \go so that \go
  evaluated at $\vl\seq{x},\seq{g}\vr$ is
\[      x_0g_{0}^{-1}dg_{0} + x_1g_{1}^{-1}dg_{1} + \, \cdots \,  +
x_ng_{n}^{-1}dg_{n},
\]
where \seq{x} are the barycentric coordinates in $\gD^n$ and
$g_{i}^{-1}dg_{i} = \pi_i^\ast(g^{-1}dg)$ for the projection $\pi_i:
G^{n+1}\ra G$ on the $i$th factor.

\item The canonical connection 1-form $E(z^{-1}dz)$ on $E\Cst$ is
  defined by the family of $E\C$-valued 1-forms $E(z^{-1}dz)^n$ on
  $\Delta^n\times (\Cst)^{n+1}$ such that  $E(z^{-1}dz)^n$ evaluated
  on a vector $(v_\gD, v_0, \dots, v_n)$ at a point  $|t_1,\, \ldots
  \, , t_n, z_0[z_1|\, \cdots \, |z_n]|$ is given by   the formula
\[\vl t_1, \, \ldots \, , t_n, z_0^{-1}v_0[z_1^{-1}v_1| \, \cdots \,
|z_n^{-1}v_n]\vr.\]
In the sequel we will use the notation
\[E(z^{-1}dz)^n_{|t_1, \, \ldots \, , t_n, z_0[z_1|\, \cdots \,  |z_n]|}=
\vl t_1, \, \ldots \, , t_n, z_0^{-1}dz_0[z_1^{-1}dz_1| \, \cdots \,
|z_n^{-1}dz_n]\vr.
\]
Similarly, the canonical connection 1-form $B(z^{-1}dz)$ on $B\Cst$ is
defined by the family of $B\C$-valued 1-forms $B(z^{-1}dz)^n$ on
$\Delta^n\times (\Cst)^n$, where
\[B(z^{-1}dz)^n_{|t_1, \, \ldots \, , t_n,
  [z_1|\, \cdots \, |z_n]|}= \vl t_1, \, \ldots \, , t_n,
[z_1^{-1}dz_1| \, \cdots \, |z_n^{-1}dz_n]\vr.\] The compatibility
conditions \eqref{compcond1}, \eqref{compcond2} are easy to check
calculations. It is also easy to see that $E(z^{-1}dz)$ is a
$E\Cst$-equivariant 1-form and $B(z^{-1}dz)$ is a $B\Cst$-equivariant
1-form.

\end{enumerate}

Let $G$ be an abelian Lie group. A smooth $p$-form on $EB^sG$ and
$B^{s+1}G$, for $s\geq 1$, is defined by the following inductive
procedure.

Suppose we have defined smooth $p$-forms on $B^{s}G$ as well as on
each product $\gD^k\times (B^{s}G)^m$ for $k,m\geq 0$. Then a smooth
$p$-form \ga on $EB^sG$ consists of a family of $p$-forms $\ga^n$ on
$\gD^n\times (B^{s}G)^{n+1}$ satisfying the compatibility conditions
\eqref{compcond1} and \eqref{compcond2}. Similarly, a smooth $p$-form
\ga on $B^{s+1}G$ consists of a family of $p$-forms $\ga^n$ on
$\gD^n\times (B^{s}G)^{n+1}$ satisfying the compatibility conditions
\eqref{compcond1} and \eqref{compcond2}. A smooth $p$-form \ga on
$\gD^k\times (B^{s+1}G)^m$ consists of a family of $p$-forms $\ga^n$
on $\gD^k\times(\gD^n\times (B^{s}G)^{n+1})^m$ satisfying the
compatibility conditions
\begin{gather*}
  (id_{\gD^k}\times (\del^i\times id)^m )^{\ast}\alpha^n =
  (id_{\gD^k}\times (id \times \del_i)^m )^{\ast}\alpha^{n-1}\\
  (id_{\gD^k}\times (s^i\times id)^m )^{\ast}\alpha^n =
  (id_{\gD^k}\times (id \times s_i)^m)^{\ast}\alpha^{n+1}
\end{gather*}

\begin{Exa}\label{e4.2}
  Now, for every $s> 0$ we are going to construct $EB^{s-1}\C$-valued
  differential 1-form $EB^{s-1}(z^{-1}dz)$ on $EB^{s-1}\Cst$ and
  $B^s\C$-valued differential 1-form $B^s(z^{-1}dz)$ on $B^s\Cst$.
  Note that from Lemma~\ref{BV} we know that for every $s> 0$ the
  groups $EB^{s-1}\C$ and $B^s\C$ are $\C$-vector spaces. Thus, it
  make sense to talk about $EB^{s-1}\C$ or $B^s\C$-valued differential
  forms.

  The canonical connection 1-form $EB^{s-1}(z^{-1}dz)$ on
  $EB^{s-1}\Cst$ is a 1-form on $EB^{s-1}\Cst$ so that
  $EB^{s-1}(z^{-1}dz)$ evaluated at $|t_1, \, \, \ldots \, \, , t_n,
  g_0[g_1|\, \cdots \, |g_n]|$ is given by the inductive formula
\[\vl t_1, \, \ldots \, , t_n, B^{s-1}(g_0^{-1}dg_0)[B^{s-1}(g_1^{-1}dg_1)|
\, \cdots \, |B^{s-1}(g_n^{-1}dg_n)]\vr,\] The canonical connection
1-form $B^s(z^{-1}dz)$ on $B^s\Cst$ is a 1-form on $B^s\Cst$ so that
$B^s(z^{-1}dz)$ evaluated at $|t_1, \, \ldots \, , t_n, [g_1|\, \cdots
\, |g_n]|$ is given by
\[\vl t_1, \, \ldots \, , t_n, [B^{s-1}(g_1^{-1}dg_1)|
\, \cdots \, |B^{s-1}(g_n^{-1}dg_n)]\vr, \] where $g_0, g_1, \, \ldots
\, , g_n \in B^{s-1}\Cst$ and $B^{s-1}(g_i^{-1}dg_i)$ is the canonical
connection 1-form $B^{s-1}(z^{-1}dz)$ on $B^{s-1}\Cst$ evaluated at
$g_i$.
\end{Exa}

\subsection{Connections on principal $B^s\Cst$-bundles}

A {\em connection} on a smooth principal $B^s\Cst$-bundle $E\ra M$ is
a collection $\{\omega_i\in \gG(U_i,B^s\A{1})\}_{i\in I}$ of
$B^s\C$-valued 1-forms, for some open covering $\{U_i\}_{i\in I}$ of
$M$, such that for every $i,j\in I$ so that $U_i\cap U_j \neq
\emptyset$
\[\omega_i -\omega_j =g_{ij}^{\ast}B^s(z^{-1}dz),\]
where $g_{ij}: U_i\cap U_j \ra B^s\Cst$ is a
transition function of the bundle  $E\ra M$.

Equivalently, a {\em connection} on a smooth principal
$B^s\Cst$-bundle $E\ra M$ is given by a $B^s\Cst$-equivariant global
section of the sheaf $B^s\cA^1_{E,\C}$.

The pull-back $g_{ij}^{\ast}B^s(z^{-1}dz)$ can be described explicitly
by the formula
$$g_{ij}^{\ast}B^s(z^{-1}dz) = \dlog (g_{ij}),$$
\n where $\dlog (g_{ij})$ \n is defined by the induction on $s$ as
follows. For any smooth function $f:U \ra B^s\Cst$ given locally by the formula
\begin{gather*}
f(x) = \vl t_1(x), \, \ldots \,  , t_n(x), [f_1(x)| \, \cdots \,  |f_n(x)]\vr,
\end{gather*}
\n where  $f_i(x): U \ra B^{s-1}\C$, we define $\dlog f \in
\gG(U,B^s\A{1})$ by
\begin{gather*}\label{dlog}
\dlog f(x) = \vl t_1(x), \, \ldots \,  , t_n(x), [\dlog f_1(x)| \,
\cdots \,  |\dlog f_n(x)]\vr.
\end{gather*}

It is easy to see that if $f:U \ra B^s\Cst$ and $\dlog (f) =0$, then
$f:U \ra B^s(\Cst)^{\delta}$, where $(\C^{\ast})^{\delta}$ is the group
$\Cst$ with the discreet topology.

\begin{Exa}\label{e2.2}\end{Exa}
\begin{enumerate}
\item It is easy to see that the differential 1-form $\go$
from Example~\ref{e2.1} is $G$-equivariant. Hence, it is a
connection 1-form on the universal principal $G$-bundle $EG \ra BG$.
We will call it {\em the canonical connection 1-form} of $EG \ra BG$.

For $G=\Cst$ the form $\go$ evaluated at
$|x_0, x_1, \, \ldots \, , x_n, z_0, z_1, \, \ldots \, , z_n|$
is given by the formula
$$\go_{|x_0, x_1, \, \ldots \, , x_n, z_0, z_1, \, \ldots \, , z_n|}=
\sum_{i=0}^n x_i \frac{dz_i}{z_i}.$$

Note that $\go$ is the composition $l\circ E\go$, where $E\go$ is
the canonical $E\C$-valued connection 1-form  on $E\Cst$ and  $l: E\C
\ra \C$ is the splitting of the \ses
$$0 \lra \C \lra E\C \lra B\C \lra 0$$
given by the formula
\begin{gather}\label{splitting}
l(|x_0, x_1, \, \ldots \, , x_n, z_0, z_1, \, \ldots \, , z_n|) =
\sum_{i=0}^n x_i z_i.
\end{gather}

\item The canonical connection 1-form on the universal $B^s\Cst$-bundle
$$B^s\Cst \lra EB^s\Cst \lra B^{s+1}\Cst$$
\n can be defined as the composition $l\circ EB^s(z^{-1}dz)$, where
$l:EB^s\C \ra B^s\C$ is the splitting of the \ses
$$0 \lra B^s\C \lra EB^s\C \lra B^{s+1}\C \lra 0$$
\n given by the formula \eqref{splitting}, where now $z_i \in B^s\C$.

\item From Corollary~\ref{smoothco} and the above example it follows that
every smooth principal $B^s\Cst$-bundle carries a connection.

In particular, the smooth principal $B\Cst$-bundle over $S^3$ from
Example~\ref{e1.1} can be equipped with the connection 1-form $\go$
that is the pull-back of the canonical connection 1-form $l\circ E\go$
on $E\Cst \ra B\Cst$. Following J-L. Brylinski and P. Deligne (see
\cite[Chapter~7]{bry-greenbook}) one can interpret \go as the Dirac
monopole.

\end{enumerate}

The ordinary exterior derivative $d: \cl{A}^k_{M,\mathbb{C}} \ra
\cl{A}^{k+1}_{M,\mathbb{C}}$ has an extension
\[      d:B^s\cl{A}^k_{M,\mathbb{C}} \lra
B^s\cl{A}^{k+1}_{M,\mathbb{C}}  \]
\n defined inductively by the formula
\[      d\bigl(\vl t_1(x), \, \ldots \, , t_n(x), [\alpha_1(x)|\, \cdots \,
|\alpha_n(x)]\vr\bigr)=
\vl t_1(x), \, \ldots \, , t_n(x), [d\alpha_1(x)|\, \cdots \,
|d\alpha_n(x)]\vr.       \]

Similarly, we define
$$d: EB^s\cl{A}^k_{M,\mathbb{C}} \lra EB^s\cl{A}^{k+1}_{M,\mathbb{C}}.$$

Note that $d \bigl(B^s(z^{-1}dz)\bigr) =0$, and hence,
\[d\omega_i -d\omega_j =dg_{ij}^{\ast}(B^s(z^{-1}dz)) =
g_{ij}^{\ast}(d B^s(z^{-1}dz)) =
g_{ij}^{\ast}(B^s(d(z^{-1}dz))) =0.
\]
Thus, the family $\{d\omega_i\}$ defines a global section \gO of the
sheaf $B^s\A{2}$, which is by definition the {\em curvature} of the
connection $\{\omega_i\}$.

\n {\bf Proof of Theorem B. } The exponential \ses
\begin{gather}\label{exp}
0\lra \u{\Z} \lra \u{\C}^\gd \lra \u{(\Cst)}^\gd \lra 0
\end{gather}
induces the \ses of the bar cochain complexes
\[ 0\lra C_B^\ast(\u{\Z}) \lra C_B^\ast(\u{\C}^\gd)
\lra C_B^\ast(\u{(\Cst)}^\gd) \lra 0            \]

For every abelian Lie group $G$, the group $H^p(M;G)$ is isomorphic to
the group $H^p(C_B^\ast(\u{G}))$, which in turn can be identified with the
group $L(M, B^{p-2}G)$ of isomorphism classes of smooth principal
$B^{p-2}G$-bundles over $M$ (see Remark~\ref{rem1}). Therefore, the
cohomology \les associated with the exponential \ses \eqref{exp}
induces a commutative diagram {\small
\[{\diagram
\rto & {H^{p-1}(M;\Z)} \rto\dto^-{\cong}&
{H^{p-1}(M;\C)} \rto\dto^-{\cong}&
{H^{p-1}(M;\Cst)} \rto\dto^-{\cong}&
{H^{p}(M;\Z)} \rto\dto^-{\cong}& \\
\rto & {L(M, B^{p-2}\Z)} \rto&
{L(M, B^{p-2}\C^\gd)} \rto&
{L(M, B^{p-2}(\Cst)^\gd)} \rto^-f&
{L(M, B^{p-2}\Cst)} \rto&
\enddiagram}    \]}

\n where the isomorphism $H^{p}(M;\Z) \ra L(M, B^{p-2}\Cst)$ is the
composition of isomorphisms
\[ H^{p}(M;\Z) \lra L(M, B^{p-1}\Z)  \lra L(M, B^{p-2}\Cst)     \]
\n and
\[      f : L(M, B^{p-2}(\Cst)^\gd) \lra L(M, B^{p-2}\Cst) \]
is the forgetful homomorphism induced by the homomorphism
$B^{p-2}(\Cst)^\gd \ra B^{p-2}\Cst$.

It is easy to see that the above diagram induces the following
commutative diagram
\[{\diagram
{0} \rto & {\frac{H^{p-1}(M;\C)}{H^{p-1}(M;\Z)_{TF}}} \rto\dto^-{\cong}&
{H^{p-1}(M;\Cst)} \rto\dto^-{\cong}&
{\tors H^{p}(M;\Z)} \rto\dto^-{\cong}& {0} \\
{0} \rto & {\ker (f)} \rto&
{L(M, B^{p-2}(\Cst)^\gd)} \rto^-f&
{\im (f)} \rto& {0}
\enddiagram}    \]
\n whose rows are exact sequences.

In order to finish the proof of Theorem~B, we have to show the for
every $q>p$ there is an isomorphism
\begin{gather}\label{basiciso}
\HD{p}{q} \cong H^{p-1}(M;\Cst)
\end{gather}
\n and that the group $L(M, B^{p-2}\Cst,\nabla^{\text{flat}})$ of
isomorphism classes of flat connections on smooth principal
$B^{p-2}\Cst$-bundles over $M$ is isomorphic to
$L(M, B^{p-2}(\Cst)^\gd)$.

The isomorphism \eqref{basiciso} follows from the fact that there is a
quasi-isomorphism
\begin{gather}\label{quasi-iso}
{\diagram
{\u{\Z(q)}} \rto\dto& {\A{0}} \rto^d\dto^{\alpha}&
{\A{1}}\rto^d\dto^{(2\pi\sqrt{-1})^{1-q}} &
{\dots} \rto^d& {\A{q-1}}\dto^{(2\pi\sqrt{-1})^{1-q}}\\
{0} \rto            & {\us{\C}} \rto^\dlog& {\A{1}} \rto^d& {\dots}
\rto^d& {\A{q-1}}
\enddiagram}
\end{gather}
where $\alpha(f) = \exp((2\pi\sqrt{-1})^{1-q}\cdot f)$, between the
smooth Deligne complex $\Z(q)^{\infty}_D$ and the complex
$\mathcal{A}^{<q}_{M,\mathbb{C}}(\dlog)[-1]$, where
\begin{align*}
\begin{CD}
\quad &\mathcal{A}^{<q}_{M,\mathbb{C}}(\dlog): \qquad
\us{\C}& @>\dlog>>
\A{1} @>\text{d}>>
\, \cdots \,  @>\text{d}>> \A{q-1}\\
\intertext{is the truncation of the complex}
\quad &\mathcal{A}^{\ast}_{M,\mathbb{C}}(\dlog): \qquad
\us{\C}& @>\dlog>>
\A{1} @>\text{d}>>
\, \cdots \,  @>\text{d}>> \A{q-1}@>\text{d}>> \, \cdots \,
\end{CD}
\end{align*}
\n which is an acyclic resolution of the constant sheaf $\Cst_M$. Therefore,
for every $q>p$ there are isomorphisms
$$\HD{p}{q} \cong \mathbb{H}^{p-1}(\mathcal{A}^{<q}_{M,\mathbb{C}}(\dlog))\cong
\mathbb{H}^{p-1}(\mathcal{A}^{\ast}_{M,\mathbb{C}}(\dlog))\cong
H^{p-1}(M;\Cst).$$

The isomorphism
\[      L(M, B^{p-2}\Cst,\nabla^{\text{flat}}) \cong
L(M, B^{p-2}(\Cst)^\gd) \]
\n is a consequence of the following lemma.

{\Le \label{flat1} There is a one-to-one correspondence between flat
connections on a smooth principal $B^{p-2}\Cst$-bundle $E\ra M$ and
reductions of the structure group of $E\ra M$ to $B^{p-2}(\Cst)^\gd$.}

\Proof Let $E\ra M$ be a smooth principal $B^{s}\C^{\ast}$-bundle with
a flat connection. We will show that the structure group of $E\ra M$
can be reduced to $B^{s}(\C^{\ast})^{\delta}$ or equivalently, that
$E\ra M$ has transition functions $\tilde{g}_{ij}: U_{ij} \ra
B^{s}(\C^{\ast})^{\delta}$.

Let $\mathcal{U}=\{U_i\}_{i\in I}$ be an open covering of $M$
consisting of contractible subsets of $M$ and let $\{g_{ij}: U_{ij}
\ra B^{s}\C^{\ast}\}_{i,j\in I}$ be a family of transition functions of $E\ra
M$. Suppose, $\{\omega_i\in \gG(U_i,B^s\A{1})\}_{i\in I}$ is a flat
connection on $E\ra M$. That is every $\omega_i$ is a closed form and
for every $i,j$ so that $U_i \cap U_j \neq \emptyset$
\[\omega_i - \omega_j =\dlog g_{ij}.\]

It is easy to see (using the Poincare Lemma and  the induction on $s$)
that if $d \omega_i = 0$, then there is a $B^s\C$-valued function
$f_i$ such that $df_i =\omega_i$.

For any smooth function $f:U\ra B^s\C$ we define, by the induction
on $s$, the $B^s\Cst$-valued function $\exp(f)$. If $f:U\ra
B^s\C$ is given locally by the formula
$$f(x) = \vl t_1(x), \, \ldots \, , t_n(x), [g_1(x)| \, \cdots \,
|g_n(x)]\vr,$$
\n where $g_i: U\ra B^{s-1}\C$, then
$$\exp f(x) = \vl t_1(x), \, \ldots \, , t_n(x), [\exp g_1(x)| \, \cdots \,
|\exp
g_n(x)]\vr.$$
Since $\dlog (\exp f) = df$,
\[\dlog g_{ij} = \omega_i - \omega_j = d(f_i - f_j) = \dlog
(\exp (f_i - f_j)) = -\dlog (\delta(\exp f)_{ij}).
\]
Therefore, for $\tilde{g}_{ij} = g_{ij} +\delta(\exp f)_{ij}$
\[\dlog \tilde{g}_{ij} = 0\]
\n and hence $\tilde{g}_{ij}: U_{ij}\ra B^{s}(\C^{\ast})^{\delta}$.
The family $\{\tilde{g}_{ij}\}$ gives the required transition
functions of $E\ra M$.

Now suppose, $E\ra M$ is a principal $B^{s}\C^{\ast}$-bundle with
transition functions $g_{ij}: U_{ij}\ra B^{s}(\C^{\ast})^{\delta}$. A
flat connection on $E\ra M$ is given by the family $\{\omega_i\}$ of
trivial (tautologicly equal to zero) 1-forms. Obviously, $d\omega_i=0$
and $\omega_i - \omega_j = 0 = \dlog g_{ij}$. \RBox

\section{$k$-connections on principal $B^s\Cst$-bundles}

In this section we define $k$-connections, $k=1,\ldots ,s+1$, and
scalar curvatures on smooth principal $B^s\Cst$-bundles and prove
Theorems A and C. In particular, we show that for $p\geq 2$ the group
$\HD{p}{p}$ is isomorphic to the group of equivalence classes of
smooth principal $B^{p-2}\Cst$-bundles with $k$-connections for
$k=1,\ldots ,p-1$.

By definition, a {\em 1-connection} on smooth principal
$B^s\Cst$-bundle is a connection on this bundle. To motivate a
definition of a $k$-connection for $k\geq 2$, we will first
reformulate the standard definition of a connection on smooth
principal $\Cst$-bundle.

A smooth principal $\Cst$-bundle $E\ra B$ is given either by a family
\[       \{g_{ij}: U_i\cap U_j \ra \Cst\}_{i,j\in I}    \]
of transition functions associated with an open covering
$\cl{U}=\{U_i\}_{i\in I}$ of $M$ or by a smooth map $g:M\ra B\Cst$ so
that $E\ra B = g^\ast (E\Cst \ra B\Cst)$. The map $g$ can be described
in terms of the transition functions $\{g_{ij}\}$ by the formula
\[      g(x) = \vl f_{i_0}(x), f_{i_1}(x),\ldots ,f_{i_{n}}(x),
[g_{i_0i_1}(x)|g_{i_1i_2}(x)|\; \cdots \;
|g_{i_{n-1}i_n}(x)]\vr,
\]
wher $\{f_i\}_{i\in I}$ is a partition of unity subordinated to the
covering $\cl{U}$ (see the proof of Proposition~\ref{smoothpr}).

Classically, a connection on $E\ra B$ is given by a family of 1-forms
$\{\go_i \in \gG(U_i,A^{1}_\C)\}_{i\in I}$ so that
\[ \go_i - \go_j = \dlog g_{ij},\qquad \text{on}\quad U_i\cap U_j\neq
\emptyset\]

Alternatively, in terms of the map $g:M\ra B\Cst$, a connection on
$E\ra B$ is a global section -\go of the sheaf $E\A{1}$ so that
$\pi_\ast \go = \dlog g$, where $\pi_\ast : \gG(M,E\A{1})\ra
\gG(M,B\A{1})$ is the homomorphism induced by the morphism of sheaves
$\pi : E\A{1}\ra B\A{1}$.

Indeed, if \go is given by
\[ \go(x) = \vl f'_{i_0}(x), f'_{i_1}(x),\, \ldots \, , f'_{i_n}(x),
\go_{i_0}(x), \go_{i_1}(x),\, \ldots \, , \go_{i_n}(x)\vr, \] then
\[  \pi_\ast\go(x) = \vl f'_{i_0}(x), f'_{i_1}(x),\, \ldots\, ,
f'_{i_n}(x), [\go_{i_0}(x): \go_{i_1}(x):\; \cdots \; :
\go_{i_n}(x)]\vr,
\]
and in the non-homogeneous coordinates
\[ \pi_\ast\go(x) = \vl f'_{i_0}(x), \ldots ,f'_{i_{n}}(x),
[\go_{i_1}(x)-\go_{i_0}(x)|\; \cdots \;
|\go_{i_n}(x)-\go_{i_{n-1}}(x)]\vr. \]

Thus, the condition $\pi_\ast \go = \dlog g$ is equivalent to the
system of equations
\[\begin{cases}
  f'_i = f_i\\ \go_j - \go_i = \dlog g_{ij}
\end{cases}\]
where the second equation holds for all $x\in M$ so that $f_i(x)\neq
0$ and $f_j(x)\neq 0$.

Since $\{f_i\}_{i\in I}$ is a partition of unity on $M$, the sets
$U_i=\{x\in M\, |\; f_i(x)\neq 0\}$ form an open covering of $M$ and
the family of 1-forms $\{-\go_i \in \gG(U_i,\cl{A}^{1}_\C)\}_{i\in I}$
determines a connection on a smooth principal \Cst-bundle induced by
the map $g:M\ra B\Cst$.

In other words, the group $\tilde{L}(M,\Cst,\na)$ of smooth principal
\Cst-bundles with connections over $M$ is the pull-back
\[\begin{CD}
  \tilde{L}(M,\Cst, \na) @>>> \gG(M,E\A{1})\\
  @VVV  @VV{\pi_\ast}V          \\ \cinf{B\Cst} @>\dlog>>
  \gG(M,B\A{1})
\end{CD}\]
of the projection $\pi_\ast : \gG(M,E\A{1})\ra \gG(M,B\A{1})$ by the
homomorphism
\[\dlog:\cinf{B\Cst}\lra \gG(M,B\A{1}). \]

The group $L(M, \Cst,\na)$ of the isomorphism classes of smooth
principal \Cst-bundles over $M$ is the quotient of $\tilde{L}(M,
\Cst,\na)$ by the action of $\cinf{E\Cst}$ given by the formula
\[f\cdot (g,\go) = (g+\pi_\ast(f),\go+\dlog(f)),\]
where $\pi_\ast:\cinf{E\Cst}\ra \cinf{B\Cst}$ is the homomorphism
induced by the projection $\pi:E\Cst\ra B\Cst$.

Essentially the same as above arguments show that a connection on a
smooth principal $B^s\Cst$-bundle induced by a map $g:M\ra
B^{s+1}\Cst$ is given by a global section -\go of the sheaf $EB^s\A{1}$
so that $\dlog g = \pi_\ast \go$, where
\[\pi_\ast: \gG(M,EB^s\A{1}) \lra \gG(M,B^{s+1}\A{1}).\]
Moreover, two pairs $(g,\go), (g',\go')\in \cinf{B^{s+1}\Cst}\oplus
\gG(M,EB^s\A{1})$ determine isomorphic smooth principal
$B^s\Cst$-bundles with connections if and only if there is a smooth
map $h:M\ra EB^s\Cst$ so that
\begin{gather*}
  (g,\go) = (g'+\pi_\ast h, \go'+\dlog h).
\end{gather*}

Now, we are going to define a 2-connection of the isomorphism class
$[E,\go]$ of a smooth principal $B\Cst$-bundle $E\ra M$ with a
connection \go.

Let $E\ra M$ be a smooth principal $B\Cst$-bundle induced from the
universal principal $B\Cst$-bundle $EB\Cst\ra B^2\Cst$ by a map
$g:M\ra B^2\Cst$ and let $\go\in \gG(M, EB\A{1})$ be a connection on
$E\ra M$. That is, $\pi_\ast\go = \dlog g$, where
$\pi_\ast:\gG(M,EB\A{1})\ra \gG(M,B^2\A{1})$ is the homomorphism
induced by the morphism of sheaves $\pi:EB\A{1}\ra B^2\A{1}$.

The curvature $-d\go$ of the connection -\go is a global section of the
sheaf $B\A{2}$, because the sequence
\[0\lra \gG(M,B\A{2}) \overset{i_\ast}{\lra} \gG(M,EB\A{2})
\overset{\pi_\ast}{\lra} \gG(M,B^2\A{2}) \lra 0   \]
is exact and $\pi_\ast (d\go) = d(\pi_\ast\go) = d(\dlog g) = 0$.

If $(g',\go')$ determines an isomorphic to $(E,\go)$ smooth principal
$B\Cst$-bundle with a connection, then $d\go =d(\go'+\dlog h) =
d\go'$, and hence a curvature determines a homomorphism
\[ d : L (M, B\Cst, \nabla) \lra \gG(M,B\A{2}), \]
where $L (M, B\Cst, \nabla)$ is the group of isomorphism classes of
smooth principal $B\Cst$-bundles with connections over $M$.

Consider the following pull-back diagram
\[      \begin{CD}
  \tilde{L}(M, B\Cst, \na_1, \na_2) @>>> \gG(M,E\A{2})\\
  @VVV  @VV{-\pi_\ast}V          \\ L(M, B\Cst, \nabla) @>d>>
  \gG(M,B\A{2})
        \end{CD}
\]
The group $\tilde{L}(M, B\Cst, \na_1, \na_2)$ consists of elements
$([g,\go_1],\go_2)$, where $g:M\ra B^2\Cst$ is a smooth map,
$[g,\go_1]$ is the isomorphism class of a smooth principal
$B\Cst$-bundle $g^\ast(EB\Cst \ra B^2\Cst)$ with a connection $-\go_1$,
and $\go_2$ is a global section of the sheaf $E\A{2}$ so that
$-\pi_\ast (\go_2) = d\go_1$. The equation $-\pi_\ast \go_2 = d\go_1$ is
an analogue of the connection condition $-\pi_\ast \go = \dlog g$,
therefore we will refer to $\go_2$ as a {\em 2-connection} of the
equivalence class $[g,\go_1]$ of the pair $(g,\go_1)$.

Note that there is an action
\[       \gG(M,E\A{1}) \times \tilde{L}(M, B\Cst,\na_1, \na_2)
\lra \tilde{L}(M, B\Cst,\na_1, \na_2)   \] of $\gG(M,E\A{1})$ on
$\tilde{L}(M, B\Cst,\na_1, \na_2)$ given by
\[\ga \cdot ([g,\go_1],\go_2) = ([g,\go_1-\gs(\ga)],\go_2+d\ga),\]
where \gs is the composition
\[ \gG(M,E\A{1}) \overset{\pi_\ast}{\lra}  \gG(M,B\A{1}) \overset{i_\ast}{\lra}
\gG(M,EB\A{1}). \] The quotient
\[L(M, B\Cst, \na_1, \na_2) = \dfrac{\tilde{L}(M, B\Cst, \na_1,
  \na_2)}{\gG(M,E\A{1})} \] will be called the group of equivalence
classes of smooth principal $B\Cst$-bundles with 1 and 2-connections
over $M$.

For $s\geq 1$, a group
\[L(M, B^s\Cst,\na_1,\na_2,\ldots , \na_{s+1}) =
L(M, B^s\Cst,\{\na_i\}_{i=1}^{s+1})\] of equivalence classes of smooth
principal $B^s\Cst$-bundles with k-connections, $k=1,2,\ldots, s+1$,
over $M$ will be defined by the following inductive procedure.

Suppose, we have already constructed the group $L(M, B^s\Cst, \{
\na_i\}_{i=1}^k)$ of equivalence classes $[g,\go_1, \ldots ,\go_k]$ of
smooth principal $B^s\Cst$-bundles with $j$-connections, for $1\leq
j\leq k<s+1$, over $M$. Then the group $L(M, B^s\Cst,
\{\na_i\}_{i=1}^{k+1})$ is defined as follows.

Consider the pull-back diagram
\[      \begin{CD}
  \tilde{L}(M, B^s\Cst, \na_1, \ldots ,\na_{k+1}) @>>>
  \gG(M,EB^{s-k}\A{k+1})\\ @VVV @VV{(-1)^k\pi_\ast}V          \\
L (M, B^s\Cst, \na_1,\ldots , \na_k) @>d>> \gG(M,B^{s-k+1}\A{k+1}),
        \end{CD}
\]
where $d([g,\go_1, \ldots ,\go_k])=d\go_k$. There is an action
\[       \gG(M,EB^{s-k}\A{k}) \times
\tilde{L}(M, B^s\Cst,\na_1, \ldots ,\na_{k+1}) \lra \tilde{L}(M,
B^s\Cst,\na_1, \ldots ,\na_{k+1})
\]
of $\gG(M,EB^{s-k}\A{k})$ on $\tilde{L}(M, B^s\Cst,\na_1, \ldots
,\na_{k+1})$ given by
\[      \ga \cdot ([g,\go_1,\ldots ,\go_k],\go_{k+1}) =
([g,\go_1,\ldots ,\go_{k-1}, \go_k+(-1)^k\gs(\ga)],\go_{k+1}+d\ga),   \]
where \gs is the composition
\[\gG(M,EB^{s-k}\A{k}) \overset{\pi_\ast}{\lra}
\gG(M,B^{s-k+1}\A{k}) \overset{i_\ast}{\lra}
\gG(M,EB^{s-k+1}\A{k}). \]
We set
\[      L(M, B^s\Cst, \na_1, \ldots ,\na_{k+1}) = \dfrac{\tilde{L}(M, B^s\Cst,
  \na_1, \ldots ,\na_{k+1})}{ \gG(M,EB^{s-k}\A{k})}.    \]

The form $(-1)^{k+1}\omega_{k+1}$, where $\go_{k+1}$ is the component
of an element $([g,\go_1,\ldots ,\go_k],\go_{k+1})$ of $\tilde{L}(M,
B^s\Cst,\{\na_i\}_{i=1}^{k+1})$ is called a {\em $(k+1)$-connection}
of $[g,\go_1,\ldots ,\go_k]$. The image of $([g,\go_1,\ldots
,\go_k],\go_{k+1})$ in $L(M, B^s\Cst,\{\na_i\}_{i=1}^{k+1})$ will be
denoted by $[g,\go_1,\ldots ,\go_{k+1}]$.

Iterating the above procedure we get the group $L(M,
B^s\Cst,\{\na_i\}_{i=1}^{s+1})$ of equivalence classes of smooth
principal $B^s\Cst$-bundles with $k$-connections, $k=1, \ldots ,s+1$.

\begin{Prop}\label{prop3.2}
  For every $p\geq 2$ there is an isomorphism
\[      \HD{p}{p} \quad \cong \quad L(M, B^{p-2}\Cst, \{\na_i\}_{i=1}^{p-1})
\]
\end{Prop}

\Proof Consider the double complex

\[ \begin{CD}
  @AAdA @AAdA   @AAdA\\ \gG(M,E\A{2}) @>>\gs> \gG(M,EB\A{2}) @>>\gs>
  \gG(M,EB^2\A{2}) @>>\gs> \\ @AAdA     @AAdA   @AAdA\\ \gG(M,E\A{1})
  @>>\gs> \gG(M,EB\A{1}) @>>\gs> \gG(M,EB^2\A{1}) @>>\gs> \\
  @AA{\dlog}A   @AA{\dlog}A     @AA{\dlog}A\\ \cinf{E\Cst} @>>\gs>
  \cinf{EB\Cst} @>>\gs> \cinf{EB^2\Cst} @>>\gs>
\end{CD} \]
of the bar cochain complexes of the components of $\A{\ast}(\dlog)$.

There is a sequence of isomorphisms
\[      \HD{p}{p} \cong \bH^{p-1}(\A{<p}(\dlog))\cong
H^{p-1}(\text{Tot}^\ast(B^{\ast,<p}_M),D), \] where
$(\text{Tot}^\ast(B^{\ast,<p}_M),D)$ is the total complex of the
double complex $B^{\ast,<p}_M = \{B^{n,s}_M\}_{s<p}$ defined as
follows

\begin{gather*}
  \tot^m(B^{\ast,<p}_M) = \bigoplus\limits_{n+s=m,
    s<p}B^{n,s}_M\\[6pt] B^{n,s}_M =
\begin{cases}
  C^\infty(M, EB^n\Cst) &\text{for}\quad s=0,n\geq 0\\
  \gG(M,EB^n\A{s})      &\text{for}\quad s>0,n\geq 0\\
  0             &\text{for}\quad s<0 \thickspace \text{or}\thickspace
  n< 0
\end{cases}     \\
\intertext{and} D: \tot^m(B^{\ast,<p}_M) \lra
\tot^{m+1}(B^{\ast,<p}_M)\\ \intertext{is so that} D =
\begin{cases}
  \dlog +\gs    &\text{on}\quad B^{n,0}_M\\ d+ (-1)^s\gs
  &\text{on}\quad B^{n,s}_M \quad \text{for $s>0$.}
\end{cases}
\end{gather*}

A $(p-1)$-cocycle in $(\text{Tot}^\ast(B^{\ast,<p}_M),D)$ is a
sequence $(g,\go_1,\ldots ,\go_{p-1})$, where $g\in
\cinf{EB^{p-1}\Cst}$ and $\go_i\in \gG(M,EB^{p-i-1}\A{i})$ so that
\[\begin{cases}
  \gs(g)=0&\\ \dlog g = \gs(\go_1)&\\ d\go_i = (-1)^{i}\gs(\go_{i+1})&
  \text{for}\; 1\leq i\leq p-2
  \end{cases}
\]

The condition $\gs(g)=0$ means that $g$ is a smooth map from $M$ to
$B^{p-1}\Cst$, the condition $\dlog g = \gs(\go_1)$ means that
$-\go_1$ is a connection on the smooth principal $B^{p-2}\Cst$-bundle
over $M$ induced by $g$, and the conditions $d\go_i =
(-1)^{i}\gs(\go_{i+1})$ mean that $(-1)^{i+1}\go_{i+1}$ is a
$(i+1)$-connection of the sequence $(g,\go_1, \ldots ,\go_i)$. It is
easy to see that two cocycles $(g,\go_1,\ldots ,\go_{p-1})$ and
$(g',\go'_1,\ldots ,\go'_{p-1})$ are cohomologous in
$(\text{Tot}^\ast(B^{\ast,<p}_M),D)$ if and only if the corresponding
principal bundles with connections are equivalent in $L(M,
B^{p-2}\Cst, \{\na_i\}_{i=1}^{p-1})$. Thus we get an isomorphism
\[ H^{p-1}(\text{Tot}^\ast(B^{\ast,<p}_M),D) \cong L(M, B^{p-2}\Cst,
\{\na_i\}_{i=1}^{p-1}). \] \qed

The rest of the section is devoted to proofs of Theorems A and C.

\n {\bf Proof of Theorem C.} Let us start from a definition of a
scalar curvature.

The {\em scalar curvature} of the element $[g,\go_1, \go_2, \ldots,
\go_{p-1}]$ of $L(M, B^{p-2}\Cst, \{\na_i\}_{i=1}^{p-1})$ is the
\C-valued $p$-form $(-1)^{p-1}d\go_{p-1}$. Note, that a priori
$d\go_{p-1}$ is a a global section of the sheaf $E\A{p}$, because
$\go_{p-1} \in \gG(M,E\A{p-1})$. But $\pi_\ast (\go_{p-1}) =
d\go_{p-2}$, and hence, $\pi_\ast (d\go_{p-1}) = d(\pi_\ast \go_{p-1})
=d(d\go_{p-2}) = 0$.  Therefore, $d\go_{p-1} \in A^p_\C (M)$.
Actually, $d\go_{p-1}$ is a closed (but not necessarily exact)
\C-valued $p$-form, because locally it is exact.

Thus, a scalar curvature induces a homomorphism
\[      s: L(M, B^{p-2}\Cst, \{\na_i\}_{i=1}^{p-1}) \lra
A^p_\C (M)_{cl},        \] \n where $A^p_\C (M)_{cl}$ is the group of
\C-valued closed $p$-forms on $M$.

The form $d\go_{p-1}$ is a $p$-cocycle in $\tot^\ast(B^{\ast,\ast}_M)$
which is cohomologous to zero in this complex, because
$d\go_{p-1}=D(g,\go, \go_2, \ldots ,\go_{p-1})$.

Since $\tot^\ast(B^{\ast,\ast}_M)$ is the acyclic resolution of
$\A{\ast}(\dlog)$, which in turn is a resolution of the constant sheaf
of the group $\Cst$, the image of $d\go_{p-1}$ in
\[ H^p(\tot^\ast(B^{\ast,\ast}_M)) \cong
H^p(\A{\ast}(\dlog)) \cong H^p(M;\Cst)  \] \n is zero. Therefore,
because the diagram
\[{\diagram
  {A^p_\C (M)_{cl}} \drto\rrto&& {H^{p}(M;\Cst)}\\ &
  {H^{p}(M;\C)}\urto & \enddiagram}     \] \n commutes and the
sequence
\[  0\lra H^{p}(M;\Z)_{TF} \lra H^{p}(M;\C) \lra H^{p}(M;\Cst) \]
\n is exact, the cohomology class of $d\go_{p-1}$ in $H^p(M;\C)$
belongs to the image of $H^p(M;\Z)$ in $H^p(M;\C)$. That is
$d\go_{p-1}$ is a closed form with integral periods. Thus, we showed
that the image $\im (s)$ of the scalar curvature homomorphism
\[      s: L(M, B^{p-2}\Cst, \{\na_i\}_{i=1}^{p-1}) \lra
A^p_\C (M)_{cl} \] is contained in the group $A^p_\C (M)_0$ of
\C-valued closed $p$-forms with integral periods on $M$.

Consider the following ``scalar curvature diagram''
\[{\diagram
  {0} \rto& {L(M, B^{p-2}\Cst, \na^{\text{flat}})} \rto\dto & {L(M,
    B^{p-2}\Cst, \{\na_i\}_{i=1}^{p-1})} \rto^-s\dto& {\im(s)}
  \rto\dto^-i& {0}\\ {0} \rto& {\H{p-1}{C^{\ast}}} \rto& {\HD{p}{p}}
  \rto& {A^{p}_{\C}(M)_0} \rto& {0} \enddiagram}\] \n where the first
vertical arrow is the isomorphisms from Theorem~B the second vertical
arrow is the isomorphisms from Propositions~\ref{prop3.2}, and $i:
\im(s) \ra A^{p}_{\C}(M)_0$ is the inclusion homomorphism.

The lower row \ses of the scalar curvature diagram is obtained from
the cohomology \les
\[      0\lra H^{p-1}(M;\Cst) \lra \bH^{p-1}(\A{<p}(\dlog))
\overset{d}{\lra} A^p_\C (M)_{cl} \lra H^{p}(M;\Cst) \lra       \] \n
associated with the \ses of sheaves
\[      0\lra \Cst_M \lra \A{<p}(\dlog) \overset{d}{\lra}
(\A{p})_{cl}[-p+1] \lra 0       \]

In order to prove the exactness of the upper row of the scalar
curvature diagram one has to show that the kernel $\ker (s)$ of the
scalar curvature homomorphism $s$ coincides with the group $L(M,
B^{p-2}\Cst,\na^{\text{flat}})$ of isomorphism classes of smooth
principal $B^{p-2}\Cst$-bundles with flat connections over $M$.

If $(g,\go_1, \ldots ,\go_{p-1})$ is a $(p-1)$-cocycle in
$\tot^\ast(B^{\ast,<p}_M)$, then the condition $d\go_{p-1}=0$ holds if
and only if $(g,\go_1, \ldots ,\go_{p-1})$ is a $(p-1)$-cocycle in
$\tot^\ast(B^{\ast,\ast}_M)$. That is $(g,\go_1, \ldots ,\go_{p-1})$
represents an element of the group $H^{p-1}(M;\Cst)$. By Theorem~B the
group $H^{p-1}(M;\Cst)$ is isomorphic to $L(M, B^{p-2}\Cst,
\na^{\text{flat}})$.  Hence
\[      \ker(s) \cong H^{p-1}(M;\Cst) \cong L(M, B^{p-2}\Cst,
\na^{\text{flat}}).     \]

It is easy to see that the scalar curvature diagram commutes.
Therefore, from 5-lemma it follows that the inclusion $i: \im(s) \ra
A^{p}_{\C}(M)_0$ is an isomorphism. This finishes the proof of
Theorem~C. \qed

\n {\bf Proof of Theorem A.} First we are going to show that there is
a commutative diagram
\[ \begin{CD}
  \HD{p}{p} @>>> \H{p}{Z}\\ @VV{\cong}V @VV{\cong}V \\ L(M, B^s\Cst,
  \{\na_i\}_{i=1}^{s+1}) @>>> L(M, B^s\Cst)
\end{CD} \]
with the vertical arrows being the isomorphisms from Propositions
\ref{prop3.2} and \ref{H(M,Z)}.

For every $s\geq 1$ there is the forgetful homomorphism
\[      \gf^L : L(M, B^s\Cst, \{\na_i\}_{i=1}^{s+1}) \lra
L(M, B^s\Cst)
\]
that sends the element $[E,\go_1, \go_2, \ldots,\go_{s+1}]$ of $L(M,
B^s\Cst, \{\na_i\}_{i=1}^{s+1})$ to the isomorphism class of the
bundle $E$. The homomorphism $\gf^L$ is surjective, because every
smooth principal $B^s\Cst$-bundle carries a connection and for every
$i\geq 1$ the homomorphism
\[ \pi_\ast :  \gG(M,EB^{s-i}\A{i}) \lra \gG(M,B^{s-i+1}\A{i})    \]
is surjective.

If $(g, \go_1, \ldots ,\go_{p-1})$ is a cocycle of
$\tot^\ast(B^{\ast,<p}_M)$, then the assignment
\[      (g, \go_1, \ldots ,\go_{p-1}) \ras \{g_{ij}\},  \]
where
\[      g(x) = \vl t_{i_1}(x), t_{i_2}(x),\ldots ,t_{i_{n}}(x),
[g_{i_0i_1}(x)|g_{i_1i_2}(x)|\ldots |g_{i_{n-1}i_n}(x)]\vr      \]
induces a homomorphism
\[ \tilde{\gf}^H : \bH^{p-1}(\A{<p}(\dlog)) \lra H^1(\u{B^{p-2}\Cst})
\]
so that the diagram
\[ \begin{CD}
  \bH^{p-1}(\A{<p}(\dlog)) @>>{\tilde{\gf}^H}> H^1(\u{B^{p-2}\Cst}) \\
  @VV{\cong}V   @VV{\cong}V     \\ L(M, B^{p-2}\Cst,
  \{\na_i\}_{i=1}^{p-1}) @>>{\gf^L}> L(M, B^{p-2}\Cst)
\end{CD} \]
\n commutes.

Composing $\tilde{\gf}^H$ with the isomorphisms
\begin{gather*}
  \HD{p}{p} \lra \bH^{p-1}(\A{<p}(\dlog))\\
\intertext{and}
  H^1(\u{B^{p-2}\Cst}) \lra \H{p}{Z}
\end{gather*}
\n we get the homomorphism
\[      \gf^H : \HD{p}{p} \lra \H{p}{Z}         \]
\n so that the diagram
\[ \begin{CD}
  \HD{p}{p} @>>{\gf^H}> \H{p}{Z} \\ @VV{\cong}V @VV{\cong}V     \\
  L(M, B^{p-2}\Cst, \{\na_i\}_{i=1}^{p-1}) @>>{\gf^L}> L(M,
  B^{p-2}\Cst)
\end{CD}\]
\n commutes.

To finish the proof of Theorem~A we have to show that there is a
commutative diagram
\[{\diagram
  {0} \rto& {\ker(\gf^L)} \rto\dto^-{\cong} & {L(M, B^{p-2}\Cst,
    \{\na_i\}_{i=1}^{p-1})} \rto^-{\gf^L}\dto^-{\cong}& {L(M,
    B^{p-2}\Cst)} \rto\dto^-{\cong}& {0}\\ {0} \rto&
  {\frac{A^{p-1}_{\C}(M)}{A^{p-1}_{\C}(M)_0}} \rto& {\HD{p}{p}} \rto&
  {\H{p}{Z}} \rto& {0} \enddiagram}\] \n with exact rows and the
vertical arrows being isomorphisms.

We have already shown that the right square of the above diagram is
commutative. Exactness of the upper row is obvious. Exactness of the
lower row \ses is derived from the cohomology \les associated with the
\ses
\[ 0 \lra \A{<p}[-1] \lra \Z(p)^\infty_D \lra \u{\Z(p)} \lra 0          \]
For details the reader is referred to the proof of Theorem~1.5.3 in
\cite{bry-greenbook}.

Now we will show that there is a homomorphism $\ker(\gf^L) \lra
A^{p-1}_{\C}(M)/A^{p-1}_{\C}(M)_0$.

Suppose $(g,\go_1, \go_2, \ldots, \go_{p-1})$ represents an element
\gL of $L(M, B^{p-2}\Cst, \{\na_i\}_{i=1}^{p-1})$ which is in the
kernel of the homomorphism $\gf^L$. That is, $g$ is a smooth map from
$M$ to $B^{p-1}\Cst$ inducing a smooth principal $B^{p-2}\Cst$-bundle
isomorphic to the trivial $B^{p-2}\Cst$-bundle over $M$. Equivalently,
$g$ is homotopic to a constant map. Hence, it has a lift to a map $h$
from $M$ into $EB^{p-2}\Cst$. That is, $\pi_\ast h =g$. Therefore, the
cocycle
\[(g,\go_1, \go_2, \ldots, \go_{p-1}) = (\pi_\ast h,\go_1, \go_2, \ldots,
\go_{p-1})\] is cohomologous to a cocycle
\[ (0, \go_1-\dlog h, \go_2, \ldots ,\go_{p-1}) =
(0, \go'_1, \go_2, \ldots ,\go_{p-1}) \]

Since the rows in the double complex $B^{\ast, \ast}_M$ are exact
(everywhere except at the zero level), there is $\gb_1\in
\gG(M,EB^{p-3}\A{1})$ so that $\gs(\gb_1) = \go'_1$. Hence, the
sequence $(0, \go'_1, \go_2,\ldots, \go_{p-1})$ is cohomologous to the
sequence $(0, 0, \go_2+d\gb_1,\ldots, \go_{p-1})$.

Iterating the above process we get a representative of \gL which is of
the form $(0,0 ,0 ,\ldots,0 ,\go_{p-1}')$. Since $\pi_\ast
(\go_{p-1}')=0$, $\go_{p-1}'$ is actually a \C-valued $(p-1)$-form on
$M$.

If $(0,0 ,0 ,\ldots,0 ,\go_{p-1}'')$ is another representative of \gL,
then there is $(\seq[p-2]{\gb})\in \tot^{p-2}(B^{\ast, \ast}_M)$ so
that
\[ (0, 0 ,0 ,\ldots,0 ,\go_{p-1}')-(0, 0 ,0 ,\ldots,0 ,\go_{p-1}'') =
D(\seq[p-2]{\gb}). \] The above equality means that $(\seq[p-2]{\gb})$
is a cocycle in $\tot^{p-2}(B^{\ast, <p-1}_M)$ whose scalar curvature
is $\go_{p-1}'-\go_{p-1}''$. From Theorem~C we know that scalar
curvatures are closed forms with integral periods. Therefore, we get a
homomorphism
\begin{gather*}
  \ker(\gf^L) \lra A^{p-1}_{\C}(M)/A^{p-1}_{\C}(M)_0 \\ [0, 0 ,0
  ,\ldots,0 ,\go_{p-1}] \ras [\go_{p-1}],
\end{gather*}
\n where $[\go_{p-1}]$ is the class of the form $\go_{p-1}$ in the
quotient $A^{p-1}_{\C}(M)/A^{p-1}_{\C}(M)_0$. It is easy to see that
this homomorphism makes the right square of the diagram of
Theorem~A commutes. Hence, by 5-lemma, it is an isomorphism. \qed

\section{Holomorphic Deligne cohomology}
\label{sec:hol_case}

In this section we define holomorphic principal \bcs{s}-bundles and
holomorphic $k$-connections on them and prove Theorem~D.

A smooth map $f: X \ra B^n\Cst$ is called a \emph{holomorphic map} if
$\db f =0$, where for
\begin{gather*}
f(x) = \vl t_1(x), \, \ldots \,  , t_n(x), [f_1(x)| \, \cdots \,
|f_n(x)]\vr,
\end{gather*}
$\db f$ is defined by the analogous to $df$ inductive formula
\begin{gather*}
\db f(x) = \vl t_1(x), \, \ldots \,  , t_n(x), [\db f_1(x)| \,
\cdots \,  |\db f_n(x)]\vr.
\end{gather*}
In a similar way we define $EB^n\Cst$-valued holomorphic maps.

A smooth principal \bcs{n}-bundle is called a \emph{holomorphic
  principal \bcs{n}-bundle} if its transition functions are
holomorphic maps. It is easy to see that if $f: X \ra B^{n+1}\Cst$ is
a holomorphic map, then the induced by $f$ principal \bcs{n}-bundle
over $X$ is a holomorphic  principal \bcs{n}-bundle. There is also an
inverse to the above statement.

\begin{Prop}\label{holpr}
  For every holomorphic principal $B^s\Cst$-bundle $E \ra M$ there is
  a holomorphic map $f:M \ra B^{s+1}\Cst$ such that $E \ra M$ is the
  pull-back of the universal principal $B^s\Cst$-bundle by $f$.
\end{Prop}

The proof of Proposition~\ref{holpr} is essentially a similar as the
proof of Proposition~\ref{smoothpr}.

Let $B^n\cl{O}_X^\ast$ and $EB^n\cl{O}_X^\ast$ be the sheaves of germs
of $B^n\Cst$ and $EB^n\Cst$-valued holomorphic maps on $X$. The
composition of the \sess
\[0\lra B^n\cl{O}_X^\ast \lra EB^n\cl{O}_X^\ast \lra
B^{n+1}\cl{O}_X^\ast \lra 0\]
gives the bar resolution
\[E\Ost \lra EB\Ost \lra EB^2\Ost \lra \cdots \]
of the sheaf $\Ost$ of non-vanishing holomorphic functions on $X$.

Let \OO{r} be the sheaf of holomorphic $r$-forms on $X$ and let
$\cl{A}^{r,s}_X$ be the sheaf of smooth $(r,s)$-forms on $X$. The
sheaf $EB^n\OO{r}$ is the kernel of the sheaf morphism
\[\db : EB^n\cl{A}^{r,0}_X \lra EB^n\cl{A}^{r,1}_X,\]
which assignes to a local section
\[\vl t_1(x), \, \ldots \, , t_n(x), \alpha_1(x)\, \ldots \,
\alpha_n(x)\vr
\]
of the sheaf $EB^n\cl{A}^{r,0}_X$ the section
\[\vl t_1(x), \, \ldots \, , t_n(x), \db\alpha_1(x)\, \ldots
\,\db\alpha_n(x)\vr
\]
of the sheaf $EB^n\cl{A}^{r,1}_X$. In the same way we define the sheaf
$B^n\OO{r}$. The composition of the \sess
\[0\lra B^n\OO{r} \lra EB^n\OO{r} \lra B^{n+1}\OO{r}\lra 0\]
gives the bar resolution
\[E\OO{r} \lra EB\OO{r} \lra EB^2\OO{r} \lra \cdots \]
of the sheaf \OO{r}.

\begin{Le}\label{ebnle}
  For every $n\geq 0$ the sheaves $EB^n\Ost$ and  $EB^n\OO{r}$ are soft.
\end{Le}

The proof of Lemma~\ref{ebnle} is essentially the same as the proof of
Lemma~\ref{acyclic}.

We will denote by $L^{hol}(X, B^{r}\Cst, \{\na_i\}_{i=1}^{q-1})$
 the group of equivalence classes of holomorphic principal
\bcs{r}-bundles over $X$ with $k$-connections, for $k=1, 2,
\ldots, q-1$, which is defined by replacing everywhere in the
definition of the group of equivalence classes of smooth principal
\bcs{r}-bundles with $k$-connections, the word ``smooth'' by the
word ``holomorphic''.

\n \textbf{Proof of Theorem D.} Let $\Omega^{<q}_X(\dlog)$ be the
complex
\[ \Ost \overset{\dlog}{\lra} \OO{1} \overset{\del}{\lra} \cdots
\overset{\del}{\lra} \OO{q-1}\]
with \Ost placed in degree zero. There is a quasi-isomorphism between
$\Omega^{<q}_X(\dlog)[-1]$ and the Deligne complex $\Z(q)_D$, which is a
holomorphic analogue of the quasi-isomorphism \eqref{quasi-iso}. Thus
\[\bH^{r}(X,\Z(q)_D) \cong \bH^{r-1}(\Omega^{<q}_X(\dlog)).\]

Consider the bar resolution $\cl{B}(\Omega^{<q}_X(\dlog))$
\[{\diagram
  & & & &\\ E\OO{2}\uto^{\del}\rto^{\gs}&
  EB\OO{2}\uto^{\del}\rto^{\gs}& & &\\ E\OO{1}\uto^{\del}\rto^{\gs}&
  EB\OO{1}\uto^{\del}\rto^{\gs}& EB^2\OO{1}\uto^{\del}\rto^{\gs} & &\\
  E\Ost \uto^{\dlog}\rto^{\gs}& EB\Ost \uto^{\dlog}\rto^{\gs}&
  EB^2\Ost \uto^{\dlog}\rto^{\gs}& EB^3\Ost \uto^{\dlog}\rto^{\gs}&
  \enddiagram}
\]
of the complex $\Omega^{<q}_X(\dlog)$. Since this is an acyclic
reolution there is an isomorphism
\[\bH^{r-1}(\Omega^{<q}_X(\dlog))\cong
H^{r-1}(\tot^\ast(B^{\ast,<q}_X)),
\]
where $B^{\ast,<q}_X$ is the global sections complex associated with
$\cl{B}(\Omega^{<q}_X(\dlog))$ and $\tot^\ast(B^{\ast,<q}_X))$ is the
total complex of $B^{\ast,<q}_X$.

A $(r-1)$-cocycle in $\text{Tot}^\ast(B^{\ast,<q}_X)$ is a
sequence $(g,\go_1,\ldots ,\go_{q-1})$, where $g\in
\gG(X, EB^{r-1}\Ost)$ and $\go_i\in \gG(X,EB^{r-i-1}\OO{i})$ so that
\[\begin{cases}
  \gs(g)=0&\\ \dlog g = \gs(\go_1)&\\ \del\go_i = (-1)^i\gs(\go_{i+1})&
  \text{for}\; 1\leq i\leq r-2
  \end{cases}\]

  The condition $\gs(g)=0$ means that $g$ is a holomorphic map from
  $X$ to $B^{r-1}\Cst$, the condition $\dlog g = \gs(\go_1)$ means
  that $-\go_1$ is a connection on the smooth principal
  $B^{p-2}\Cst$-bundle over $X$ induced by $g$, and the conditions
  $\del\go_i = (-1)^i\gs(\go_{i+1})$ mean that $(-1)^{i+1}\go_{i+1}$
  is a $(i+1)$-connection of the sequence $(g,\go_1, \ldots ,\go_i)$.
  Two cocycles $(g,\go_1,\ldots ,\go_{q-1})$ and $(g',\go'_1,\ldots
  ,\go'_{q-1})$ are cohomologous in $\text{Tot}^\ast(B^{\ast,<q}_X)$
  if and only if the corresponding principal bundles with connections
  are equivalent. This gives us an isomorphism
\[ H^{r-1}(\text{Tot}^\ast(B^{\ast,<q}_X)) \cong L^{hol}(X, B^{r-2}\Cst,
\{\na_i\}_{i=1}^{q-1}). \]

In order to get the commutative diagram from Theorem~D consider the
bar resolution
\begin{gather}\label{holbarres}
0\lra \cl{B}(\Omega^{<p}_X[-1])\lra \cl{B}(\Z(p)_D) \lra
\cl{B}(\underline{\Z(p)}_X) \lra 0
\end{gather}
the \ses
\[0\lra \Omega^{<p}_X[-1] \lra \Z(p)_D \lra
\underline{\Z(p)}_X \lra 0\]

Since
\[\bH^{2p}(\tot^\ast(\cl{B}(\Z(p)_D))) \cong L^{hol}(X,
B^{r}\Cst, \{\na_i\}_{i=1}^{p-1}), \]
and
\[\bH^{2p}(\tot^\ast(\cl{B}(\underline{\Z(p)}_X))) \cong L(X,
B^{r}\Cst). \]
The hypercohomology \les associated with \eqref{holbarres} induces the
lower \ses of the diagram from Theorem~D. The quasi-isomorphisms
\begin{gather*}
\Omega^{<p}_X[-1] \lra \tot^\ast(\cl{B}(\Omega^{<p}_X[-1]))\\
  \Z(p)_D \lra \tot^\ast(\cl{B}(\Z(p)_D))\\
\underline{\Z(p)}_X \lra \tot^\ast(\cl{B}(\underline{\Z(p)}_X))
\end{gather*}
induce the vertical isomorphisms in this diagram. \qed

\setcounter{section}{1}
\setcounter{Th}{0}
\renewcommand{\thesection}{\Alph{section}}

\section*{Appendix A\\Principal Bundles, Topological Extensions, and Gerbs}
\label{sec:appA}

In this appendix we show that there is an isomorphism between the
group of isomorphism classes of smooth (or holomorphic) principal
\bcs{}-bundles over a manifold $M$ (or a complex projective variety
$X$) and the group of equivalence classes of smooth (or holomorphic)
gerbes bound by \us{\C} (or \Ost). This isomorphism is induced by a
construction, described in \cite{bry-greenbook}, which assigns to a
principal $G$-bundle $\pi: E\ra B$ and a topological central extension
\[1\lra C \lra K \lra G\lra 1\]
a sheaf of goupoids $\cl{G}_\pi$ measuring the
obstruction to the existence of a reduction of the structure group of
$\pi: E\ra B$ to $K$ (see pp.  171-172 in \cite{bry-greenbook}). In
the case of smooth (or holomorphic) principal \bcs{}-bundles and the
extension
\[0\lra \Cst \lra E\Cst \lra \bcs{}\lra 0\]
the gerbe $\cl{G}_\pi$ is equivalent to the gerbe of sections of the
bundle.  We will also describe a procedure which assigns to connection
on a principal \bcs{}-bundle a connective structure on the associated
gerbe (see pp. 169-170 in \cite{bry-greenbook}).

Let us start with a definition of a gerbe. A {\em gerbe} on a space $X$
is a sheaf of categories \cl{C} on $X$ (for the precise definition of
a sheaf of categories see Chapter~5 in \cite{bry-greenbook})
satisfying the following three conditions
\begin{itemize}
\item For every open subset $U\subset X$ the category $\cl{C}(U)$ is a
groupoid, that is, every morphism is invertible.
\item Each point $x\in X$ has a neighborhood $U_x$ for which
$\cl{C}(U_x)$ is non-empty.
\item Any two objects $P_1$ and $P_2$ of $\cl{C}(U)$ are locally
isomorphic. This means that each $x\in U$ has a  neighborhood $V$ such
that the restrictions of $P_1$ and $P_2$ to $V$ are isomorphic.
\end{itemize}

A gerbe \cl{C} is said to be {\em bound} by a sheaf \cl{A} of abelian
groups on $X$, if for every open set $U\subset X$ and every object $P$
of $\cl{C}(U)$ there is an isomorphism of sheaves
\[\alpha: \underline{\text{Aut}}(P) \lra \cl{A}|_U,\]
where $\cl{A}|_U$ is the restriction of the sheaf \cl{A} to $U$, and
$\underline{\text{Aut}}(P)$ is the sheaf of authomorphisms of $P$ so
that for an open subset $V$ of $U$ the group
$\underline{\text{Aut}}(P)(V)$ is the group of authomorphisms of the
restriction $r_V(P)$ of $P$ to $V$.  Such an isomorphism is supposed
to commute with morphisms of \cl{C} and must be compatible with
restriction to smaller open sets.

Two gerbes \cl{C} and \cl{D}  bound by \cl{A} on a manifold $M$ are
\emph{equivalent} if the following two conditions are satisfied.
\begin{itemize}
\item For every open subset $U$ of $M$ there is an equivalence of
  categories $\phi(U): \cl{C}(U) \ra \cl{D}(U)$ so that for every
  object $P$ of $\cl{C}(U)$ there is a commutative diagram
\[{\diagram
\text{Aut}_{\cl{C}(U)}(P) \rrto^-{\phi(U)} \drto_{\alpha_\cl{C}}&&
\text{Aut}_{\cl{D}(U)}(P) \dlto^{\alpha_\cl{D}} \\
& \gG(U,\cl{A}) &
\enddiagram}\]

\item For every pair of open subsets $V, U$ of $M$ so that $V \subset
  U$ there is an invertible natural transformation
\[\beta: \phi(U)\circ r_{\cl{D}} \lra r_{\cl{C}}\circ \phi(V),\]
where
\[r_{\cl{C}}: \cl{C}(U) \lra \cl{C}(V), \qquad r_{\cl{D}}: \cl{D}(U)
\lra \cl{D}(V),\] are the restriction natural transformations. It is
required that for a triple of open set $V \subset U\subset W $ in $M$
some compatibility conditions are satisfied (see p. 200 in
\cite{bry-greenbook}).
\end{itemize}

With every principal $G$-bundle $\pi: E\ra B$ and every central
extension of topological groups
\[1\lra C \lra K \lra G\lra 1\]
we can associate a gerbe $\cl{G}_\pi$ bound by \u{C} on $B$.
The gerbe $\cl{G}_\pi$ is derived from the sheaf of sections of the
bundle $\pi: E\ra B$.  For every open subset $U$ of $B$ the objects
and morphisms of $\cl{G}_\pi(U)$ are defined as follows.

Every section $s: U\ra \pi^{-1}(U)$ of $\pi^{-1}(U)\ra U$ can be
identified with a $G$-equivariant map
\[t_s: \pi^{-1}(U)\lra G\]
so that for every $\xi \in \pi^{-1}(U)$ we have $t_s(\xi)\cdot
s(\pi(\xi)) =\xi$. Let $E_s \ra \pi^{-1}(U)$ be the pull-back of
principal $C$-bundle $K\ra G$ from $G$ to $\pi^{-1}(U)$, by the map
$t_s: \pi^{-1}(U) \ra G$. It is clear that the composition
$\pi\circ\pi_s : E_s \ra U$ is a principal $K$-bundle, and hence a
reduction of the structure group of $\pi^{-1}(U)\ra U$ to $K$. The
objects of $\cl{G}_\pi(U)$ are pairs $(E,f)$ of principal $K$-bundles
$\tilde{\pi}: E\ra U$ and principal $C$-bundles $f: E \ra \pi^{-1}(U)$
so that the diagram
\[{\diagram
E \rrto^-f \drto_{\tilde{\pi}}&&
\pi^{-1}(U) \dlto^\pi \\
& U &
\enddiagram}\]
\n commutes. A morphism from $(E,f)$ to $(E',f')$ is a morphism of
principal $K$-bundles $g: E\ra E'$ so that the diagram
\[{\diagram
E \rrto^g\drto_f&& E'\dlto^{f'}\\
& \pi^{-1}(U) &
\enddiagram}
\]
commutes. The above condition implies that the group of authomorphisms
of any object $(E,f)$ of $\cl{G}_\pi(U)$ is the group of maps from $U$
to $C$, which is the section of the sheaf \u{C} over $U$. Thus
$\cl{G}_\pi$ is the gerbe bound by \u{C}.

Note, that the gerbe $\cl{G}_\pi$ has a global section if and only if
there is a reduction of the structure group of $\pi: E\ra B$ to $K$.
In particular, if $\pi: E\ra B$ is a principal \bcs{}-bundle, and our
extension is the universal extension
\[0\lra \Cst \lra E\Cst \lra \bcs{}\lra 0\]
then the associated with $\pi: E\ra B$ gerbe $\cl{G}_\pi$ measures the
obstruction for the existence of a reduction of the structure group of
$\pi: E\ra B$ to $E\Cst$. Since $E\Cst$ is contractible, every
principal $E\Cst$-bundle is trivial. Thus, the gerbe  $\cl{G}_\pi$ has a
global section if and only if $\pi: E\ra B$ is a trivial \bcs{}-bundle.
The same property has the gerbe $\cl{S}_\pi$ of local sections of the
bundle $\pi: E\ra B$, which is defined as follows. For every open
subset $U$ of $M$ the objects of $\cl{S}_\pi(U)$ are sections of $\pi:
E\ra B$ over $U$. Every local section $s: U\ra\pi^{-1}(U)$ of $\pi:
E\ra B$ induces a \bcs{}-equivariant map $t_s: \pi^{-1}(U)\ra
\bcs{}$, which in turn gives a map $\tau_s = s\circ t_s: U \ra
\bcs{}$. Let $L_s$ be the principal \Cst-bundle over $U$ induced by
the map $\tau_s$. A morphism between the objects $s,s'\in
\cl{S}_\pi(U)$ is a morphism $L_s\ra L_{s'}$ of the corresponding
principal \Cst-bundles. It is clear that $\cl{S}_\pi$ is a gerbe bound
by \us{\C}. It is not a difficult exercise to see that the natural
transformation $\cl{S}_\pi(U)\ra \cl{G}_\pi(U)$ sending a section $s$
to the pull-back $E_s$ of the universal principal \Cst-bundle by $t_s$
is an equivalence of categories that extends to an equivalence of
gerbes $\cl{S}_\pi\ra \cl{G}_\pi$.

The following theorem is an easy consequence of Theorem~H (see the
introduction) and Theorem~5.2.8 from \cite{bry-greenbook}.

\begin{Th}
  A map which sends to the isomorphism class of a principal
  $B\Cst$-bundle $\pi: E\ra B$ the equivalence class of the gerbe of
  section $\cl{S}_\pi$ of $\pi: E\ra B$ induces an isomorphism between
  the group of isomorphism classes of principal $B\Cst$-bundles and
  the group of equivalence classes of gerbes bound by \Cst.
\end{Th}

Let \cl{G} be a gerbe on $M$ bound by \us{\C}. A {\em connective
structure} on \cl{G} is an assignment to each object $P$ in
$\cl{G}(U)$ a \A{1}-torsor $\cl{Co}_P$ on $U$. That is
$\cl{Co}_P$ is a sheaf with an action of \A{1} on $\cl{Co}_P$ such that
every point has a \nbhd $U$ with the property that for each open set
$V\subset U$ the group $\cl{Co}_P(V)$ is a principal homogeneous space
under the group $\gG(V,\A{1})$. The assignment $P \ras \cl{Co}_P(U)$
should be functorial with respect to restriction of $U$ to smaller
open set and should be so that for any morphism $\gp:P\ra Q$ of objects
of $\cl{G}(U)$ (necessarily an isomorphism since \cl{G} is a gerbe),
there is an isomorphism $\gp_\ast: \cl{Co}_P(U) \ra \cl{Co}_Q(U)$ of
\A{1}-torsors, which is compatible with composition of morphisms in
$\cl{G}(U)$ and also compatible with restrictions to smaller open
sets. If \gp is an automorphism of $P$ induced by a \Cst-valued
function $g$, we require that $\gp_\ast$ be the automorphism $\nabla
\ras \nabla -\frac{dg}{g}$ of the \A{1}-torsor $\cl{Co}_P(U)$.  In a
similar way one can define a holomorphic connective structure on a
holomorphic gerbe bound by \Ost.

A connection $\go$ on a smooth principal $B\Cst$-bundle $\pi:E\ra M$
induces the following connective structure on $\cl{G}_\pi$. Let $U$ be
an open subset of $M$ so that $\cl{G}_\pi (U)$ is non-empty and let
$\go_U$ be the restriction of $\go$ to $\pi^{-1}(U)$. To every element
$(E,f)$ of $\cl{G}_\pi (U)$ we assign a set $\cl{Co}_E^\omega(U)$ of
connections on $E$ compatible with \go. That is $\tilde{\go} \in
\cl{Co}_E^\omega(U)$ if $q\circ \go = f^\ast \go$, where $q: E\C \ra
B\C$ and $f: E \ra \pi^{-1}(U)$ is the principal \Cst-bundle. It is
easy to see that the assignment $\omega \ras \cl{Co}^\omega$ is a
connective structure on $\cl{G}_\pi$ (for detail see pp. 169-170 in
\cite{bry-greenbook}). The equivalence of gerbes $\cl{S}_\pi\ra
\cl{G}_\pi$ can be used to pull-back the connective structure from
$\cl{G}_\pi$ to $\cl{S}_\pi$.  A similar to the above construction
assigns to a holomorphic connection on a holomorphic principal
$B\Cst$-bundle $E\ra X$ a holomorphic connective structure on the
associated with $E\ra X$ holomorphic gerbe.

\begin{Th}
  A map which sends to the isomorphism class of a principal
  $B\Cst$-bundle $\pi: E\ra B$ with a connection $\omega$ the
  equivalence class of the gerbe of section $\cl{S}_\pi$ of $\pi: E\ra
  B$ with the connective structure on $\cl{S}_\pi$ induced by $\omega$
  induces an isomorphism between the group of isomorphism classes of
  principal $B\Cst$-bundles with connection and the group of
  equivalence classes of gerbes bound by \Cst with connective
  structures.
\end{Th}

We leave the proof of this theorem as an exercise for the reader.

\setcounter{section}{2}
\section*{Appendix B\\The Geometric Bar Construction}

The objective of this appendix is twofold. First, we define and review
basic properties of the geometric bar construction. Second, we explain
how the geometric bar construction can be derived from the projective
space construction. Our basic references for the geometric bar
construction are \cite{milg-bar} and the survey paper \cite{sta-hsp}.

The geometric bar construction assigns to every topological group $G$
a sequence of principal $G$-bundles $E_n\ra B_n$
\[{\diagram
  {G} \dto\rdouble& {E_1}\dto\rh&{E_2}\dto\rh&{\; \cdots\;}
  \rh&{E_n}\dto\rh&{\; \cdots} \\ {pt} \rdouble& {B_1}\rh& {B_2}
  \rh&{\; \cdots \;} \rh&{B_n} \rh&{\; \cdots} \enddiagram}\] so that
for every $n\geq 0$ the space $E_n$ is contractible in $E_{n+1}$. The
universal principal $G$-bundle $EG \ra BG$ is the union
$\bigcup_{n\geq 1} E_n \ra \bigcup_{n\geq 1} B_n$ taken with the weak
topology.

If $G$ is an abelian topological group, then $EG$ and $BG$ are abelian
topological groups and the projection $EG\ra BG$ is a continuous
homomorphism with $G$ as the kernel.

The geometric bar construction is functorial and it preserves
products. That is, every continuous homomorphism $f:G\ra H$ induces
continuous maps $Ef: EG\ra EH$ and $Bf: BG\ra BH$, which are
homomorphisms for $G$ abelian, and
\[E(G\times H) = EG\times EH\qquad B(G\times H) = BG\times BH,\]
where each product is taken with the compactly generated
topology\footnote{To every topological space $X$ one can assign a
  space $(X,k)$ with compactly generated topology so that a set is
  open in $(X,k)$ if and only if its intersection with every compact
  subset of $X$ is open.}.

If $G$ is a countable CW-group\footnote{A topological group $G$ is
  called a countable CW-group if it is a countable CW-complex so that
  the map $g\ras g^{-1}$ of $G$ into itself and the product map
  $G\times G \ra G$ are both cellular (that is, they carry the
  $k$-skeleton into the $k$-skeleton).}, then $EG$ and $BG$ are
countable CW-complexes. In this appendix $G$ is a countable CW-group.
Actually, in the main body of the paper $G$ is \Z, \C, \Cst, $S^1$, or
the abelian group of a separable \C-vector space. A techincal
advantage of working with countable CW-groups is that on the spaces
appearing in the definitions of $EG$ and $BG$ one can take the
product, versus compactly generated, topology.

The archetypes of the geometric bar construction are infinite real,
complex, and quaternionic projective spaces. Actually, Milnor found a
construction that associates with every topological group $G$ a
principal $G$-bundle $E_\gD G \ra B_\gD G$, which is a limit of a
sequence of principal $G$-bundles $(E_\gD G)_n \ra (B_\gD G)_n$, so
that for $G=S^0, S^1$, and $S^3$ the bundle $(E_\gD G)_n \ra (B_\gD
G)_n$ is isomorphic to $S^n\ra \R\P^n, S^{2n+1}\ra \C\P^n$, and
$S^{4n+3}\ra \bH\P^n$ respectively.

A drawback of Milnor's constraction is that for $G$ being an abelian
topological group the spaces $E_\gD G$ and $B_\gD G$ are not abelian
groups, so the construction cannot be iterated. The geometric bar
construction is a ``normalized version'' of Milnor's construction that
fixes this problem.

There are several approaches to geometric bar construction (for a
survey on this subject see \cite{sta-hsp}).  Usually, the spaces $EG$
and $BG$ are defined as the quotients of the disjoint unions
$\coprod_{n\geq 0} \gD^n\times G^{n+1}$ and $\coprod_{n\geq
  0}\gD^n\times G^{n}$ respectively, by cetrain equivalence relations.
To explain the geometric meaning of these relations we preceded the
formal definition of geometric bar construction with the Milnor and
the Dold-Lashof constructions \cite{miln-cub2}, \cite{dol&las-pqf}.

\subsection{The unnormalized geometric bar construction.}

Let $G$ be a countable CW-group. The {\em join} $G\ast G$ is the
quotient of the product $\gD^1\times G\times G$ of the standard
1-simplex
\[\gD^1 = \{(x_0,x_1)\in \R^2 |\; x_0,x_1\geq 0,\quad
x_0+x_1=1\}\] with $G\times G$ by the equivalence relation
\[(0,1,g_0,g_1)\sim (0,1,e,g_1),\qquad (1,0,g_0,g_1) \sim (1,0,g_0,e)
\]
where $e$ is the neutral element of $G$. The equivalence class of the
sequence $(x_0,x_1,g_0,g_1)$ will be denoted by $x_0g_0\oplus x_1g_1$.

Let $I=[0,1]$. The homeomorphism
\[ \gD^1\times G\times G\lra G\times I\times G,\quad
(x_0,x_1,g_0,g_1)\ras (g_0, x_1, g_1) \] induces a homeomorphisms
between $G\ast G$ and the quotient of the product $G\times C(G)$ of
$G$ with the cone $C(G)=I\times G/0\times G$ by the equivalence
relation
\[      (g_0,[1,g_1]) \sim (e,[1,g_1])  \]
where $[t,g]$ is the image of the pair $(t,g)$ in $C(G)$.

The $(n+1)$-fold join
\[      G\ast (n+1) \ast G = \underbrace{G\ast (G\ast \, \cdots \,  (G\ast
  G)}_{(n+1)\; \text{times}} \, \cdots \, ) \] which we will also
denote by $(E_\gD G)_n$, can be identified with the quotient of the
product $\gD^n\times G^{n+1}$ of the standard simplex
\[      \gD^n = \{ (\seq{x})\in \R^{n+1} |\; x_i\geq 0, \quad
\sum\limits_{i=0}^n x_i =1 \} \] and the $(n+1)$-fold product
$G^{n+1}$ of $G$ with itself, by the equivalence relation
\begin{align*}
  &(\seq[i-1]{x}, 0, x_{i+1}, \, \ldots \, , x_n, \seq[i]{g}, \,
  \ldots \, , g_n) \sim \\ \sim &(\seq[i-1]{x}, 0, x_{i+1}, \, \ldots
  \, , x_n, g_0, \, \ldots \, ,e, \, \ldots \, , g_n)
\end{align*}

Actually, if we denote by $x_0g_0\oplus \, \cdots \, \oplus x_ng_n$
the equivalence class of the sequence $(\seq{x},\seq{g})$, then the
homeomorphism between $(E_\gD G)_n$ and the quotient $(\gD^n\times
G^{n+1})/\sim$ is given by
\begin{gather*}
  x_0g_0\oplus (1-x_0)\Bigl(x_1g_1\oplus (1-x_1)\bigl( \, \cdots \,
  (x_{n-1}g_{n-1}\oplus (1-x_{n-1})g_n) \, \cdots \,
  \bigr)\Bigr)\ras\\ \ras x_0g_0\oplus (1-x_0)x_1g_1\oplus \, \cdots
  \, \oplus \prod\limits_{i=0}^{n-1}(1-x_i)g_n
\end{gather*}

On the other hand, $(E_\gD G)_n = G\ast (E_\gD G)_{n-1}$ can be
identified with the quotient of $G\times C((E_\gD G)_{n-1})$ by the
quivalence relation
\[      (g,[1,x]) \sim (e,[1,x]).       \]
Note that there are inclusions
\[{\diagram
  {(E_\gD G)_{n-1}} \rh_-i& {C((E_\gD G)_{n-1})} \rh_-j& {(E_\gD
    G)_{n}} \enddiagram}\] given by $i(y)=[1,y]$ and $j([t,y]) =
(1-t)e\oplus ty$.

Let
\[ E_\gD G = \bigcup_{n\geq 2} (E_\gD G)_n = \bigcup_{n\geq 2}
C((E_\gD G)_n). \]

Since $E_\gD G$ is the union of cones, it is contractible.

There is a free action of $G$ on $(E_\gD G)_n$ given by
\begin{gather}\label{action}
  g\cdot (x_0g_0\oplus \, \cdots \, \oplus x_ng_n) = x_0(g g_0)\oplus
  \, \cdots \, \oplus x_n(g g_n)
\end{gather}
The orbit space of this action is denoted by $(B_\gD G)_n$. For
example, $(B_\gD G)_0$ is a single point and $(B_\gD G)_1$ is the
suspension of $G$.

Since the actions of $G$ on $(E_\gD G)_n$ and $(E_\gD G)_{n+1}$ are
compatible with the embedding $(E_\gD G)_n \subset (E_\gD G)_{n+1}$,
there is a free action of $G$ on $E_\gD G$. The quotient space $(E_\gD
G)/G$ is denoted by $B_\gD G$ and the natural map $E_\gD G \ra B_\gD
G$ is Milnor's universal principal $G$-bundle.

\begin{Exa}\label{appex1}
  For $G=S^0 = \Z/2\Z=\{\pm 1\}$ there are homeomorphisms
\[(E_\gD S^0)_n \cong S^n,\qquad (B_\gD S^0)_n \cong \R\P^n     \]
induced by the map
\[(E_\gD S^0)_n\ni x_0\ga_0\oplus \, \cdots \,  \oplus x_n\ga_n \ras
(\ga_0\sqrt{x_0}, \, \ldots \, , \ga_n\sqrt{x_n}) \in S^n \]
Similarly, for $G=S^1=U(1), S^3=SU(2)$, or \Cst there are the
following homeomorphisms
\begin{align*}
  (E_\gD S^1)_n &\cong S^{2n+1},& (B_\gD S^1)_n &\cong\C\P^n\\ (E_\gD
  S^3)_n &\cong S^{4n+3},& (B_\gD S^3)_n &\cong\bH\P^n\\ (E_\gD
  \Cst)_n &\cong S^{2n+1}\times (\R_{+})^{n+1},& (B_\gD\Cst)_n
  &\cong\C\P^n\times S^n_{+}
\end{align*}
where $\R_{+}$ is the set of positive real numbers and $S^n_{+}$ is
the intersection $(\R_{+})^{n+1}\cap S^n$.
\end{Exa}

Sometimes, it is convenient to replace the diagonal action
\eqref{action} of $G$ on $(E_\gD G)_n$ by the action of $G$ on the
first component of $(E_\gD G)_n$. This can be done by introducing a
{\em non-homogeneous coordinates} on $(E_\gD G)_n$
\[ \vl \seq{x},h_0[h_1| \, \cdots \,  |h_n]\vr_\gD =
x_0g_0\oplus x_1(g_0^{-1}g_1)\oplus \, \cdots \, \oplus x_n
(g_{n-1}^{-1}g_n)\] where $(\seq{g})\in G^{n+1}$, $h_0 =g_0$, and $h_i
= g_{i-1}^{-1}g_i$ for $i>0$.

The non-homogeneous coordinates on $(E_\gD G)_n$ lead to yet another
model of $(E_\gD G)_n$, due to Dold and Lashof \cite{dol&las-pqf}.
For example, in the 2-fold join $G\ast G$ the relations
\[ 0g_0\oplus 1g_1 = 0e\oplus 1g_1,\qquad 1g_0\oplus 0g_1 = 1g_0\oplus
0e \] correspond, in the non-homogeneous coordinates, to the relations
\[ \vl 0,1,h_0[h_1]\vr_\gD = \vl 0,1,e[h_0h_1]\vr_\gD,\quad
\vl 1,0, h_0[h_1]\vr_\gD = \vl 1,0, h_0[e]\vr_\gD \]

Thus, the symbols $\vl x_0, x_1 , h_0[h_1]\vr_\gD$ can be identified
with the points of the space
\[DL(G)=G\times C(G)\cup_\mu G =
(G\times C(G)\sqcup G)/\sim \] where $\mu:G\times G \ra G$ is the
group operation in $G$ and $\sim$ is an equivalence relation
identifying $(h_0,[1,h_1])$ with $\mu(h_0,h_1)=h_0h_1$.

Note, that there is an action of $G$ on $DL(G)$ given by
\[g\cdot (\vl x_0, x_1 , h_0[h_1]\vr_\gD) = \vl x_0, x_1 ,
gh_0[h_1]\vr_\gD\] and hence we can apply the above construction to
$DL(G)$.

In general, to any space $E$ with a $G$-action $\mu:G\times E\ra E$ we
can associate the space
\[ DL(E) =  G\times C(E)\cup_\mu G = (G\times C(E)\sqcup G)/\sim \]
where $\sim$ is an equivalence relation identifying $(h,[1,x])$ with
$\mu(h,x)$, and the action
\[ G\times DL(E)\ra DL(E), \qquad g\cdot (h|t|x) = (gh)|t|x     \]
where $h|t|x$ is the equivalence class of the sequence $(h,[t,x])\in
G\times C(E)$ in $DL(E)$.

The spaces $DL(E)$ and $G\ast E$ are $G$-equivariantly homeomorphic to
each other with the $G$-equivariantly homeomorphisms given by
\begin{align*}
  &DL(E) \ra G\ast E, & h|t|y &\ras th \oplus (1-t)(hy)\\ &G\ast E\ra
  DL(E), & x_0h\oplus x_1y &\ras h|x_0|h^{-1}y
\end{align*}
Therefore, the bundles $DL(E)\ra DL(B)$ and $G\ast E \ra (G\ast E)/G$
are isomorphic.

Applying $n$ times the Dold-Lashof construction to a topological group
$G$, we get a principal $G$-bundle $DL^n(G)\ra DL^n(G)/G$ which is
isomorphic to the bundle $(E_\gD G)_n \ra (B_\gD G)_n$.

\subsection{The geometric bar construction.} In general, the spaces
$E_\gD G$ and $B_\gD G$ have not group structure, but for $G$
abelian, some quotients of these spaces are groups.

The appropriate quotients are obtained by replacing the cone $C(E)$ in
the Dold-Lashof construction $DL(E)$ by the reduced cone
\[ \tC(E) = (I\times E)/(0\times E \cup I\times e)      \]
where $e$ is a base point of $E$. For example, for $(E,e) = (G,e)$
where $G$ is a topological group with the neutral element $e$ we
define
\[ \wt{DL}(G) = G\times \tC(G)\cup_\mu G        \]
The space $\wt{DL}(G)$ is a quotient of $DL(G)$ by the equivalence
relation
\[ h|t|e = h|0|e.       \]
The  group action of $G$ on $DL(G)$ decents to a group action of $G$
on $\wt{DL}(G)$. Thus, we can iterate this construction getting for
every $n\geq 1$ a space $\wt{DL}^n(G)$ with a free action of $G$ on
itself. We set $(EG)_n = \wt{DL}^n(G)$ and $(BG)_n = \wt{DL}^n(G)/G$.

It is easy to see that $(EG)_n$ is the quotient of the disjoint union
$\coprod\limits_{m=0}^n \gD^m\times G^{m+1}$
by the equivalence relations
\begin{multline*}
(\seq[m]{x}, \seq[m]{g}) \sim\\
\sim \begin{cases}
(\seq[i]{x}+x_{i+1}, \, \ldots \,  ,x_m, g_0, \, \ldots \,  ,\hat{g_i}, \,
\ldots \,
,g_m)&\text{for $g_i=g_{i+1}$ or $x_i=0, \; 0\leq i<m$}\\
(\seq[m-1]{x}+x_{m},\seq[m-1]{g})&\text{for $g_{m-1}=g_{m}$ or
$x_m=0$}
\end{cases}
\end{multline*}

In the non-homogeneous coordinates on $(EG)_n$ the above relations
take the form
\begin{multline*}
(t_1, \, \ldots \,  , t_m, h_0[h_1| \, \cdots \,  |h_m]) \sim\\
\sim \begin{cases}
(t_2, \, \ldots \,  , t_m, h_0h_1[h_2| \, \cdots \,  |h_m])&\text{for $t_1=0$
or
$h_0=e$} \\
(t_1, \, \ldots \,  , \hat{t_i}, \dots , t_m, h_0[h_1| \, \cdots \,
|h_ih_{i+1}|
\, \cdots \,  |h_m])&
\text{for $t_{i}=t_{i+1}$ or $h_i=e$} \\
(t_1, \, \ldots \,  , t_{m-1}, h_0[h_1| \, \cdots \,  |h_{m-1}])& \text{for
$t_m=1$
or $h_m=e$}
\end{cases}
\end{multline*}
where $0\leq t_1 \leq t_2 \leq \, \cdots \,  \leq t_m \leq 1$ are
non-homogeneous coordinates on $\gD^n$ related with the baricentric
coordinated $\seq{x}$ on $\gD^n$  by the formula
\[ t_i = x_0+x_1+ \, \cdots \,  + x_{i-1}.      \]

The equivalence class of a sequence $(\seq[m]{x}, \seq[m]{g})$ will be
denoted by $\vl \seq[m]{x}, \seq[m]{g} \vr$ and the equivalence class
of a sequence $(t_1, \, \ldots \,  , t_m, h_0[h_1| \, \cdots \,  |h_m])$ will
be
denoted by $\vl t_1, \, \ldots \,  , t_m, h_0[h_1| \, \cdots \,  |h_m]\vr$.

The space $EG$ is the quotient of the disjoint union
$\coprod\limits_{m=0}^\infty \gD^m\times G^{m+1}$
by the above equivalence relations.

Similarly, $(BG)_n$ is the quotient of the disjoint union
$\coprod\limits_{m=0}^n \gD^m\times G^{m}$
by the equivalence relations
\begin{multline*}
(\seq[m]{x}, [g_0: \, \cdots \,  :g_m]) \sim\\
\sim \begin{cases}
(\seq[i]{x}+x_{i+1}, \, \ldots \,  ,x_m, [g_0: \, \cdots \,  :\hat{g_i}: \,
\cdots \,
:g_m])&\text{for $g_i=g_{i+1}$ or $x_i=0, 0\leq i<m$}\\
(\seq[m-1]{x}+x_{m},[g_0: \, \cdots \,  :g_{m-1}])&\text{for $g_{m-1}=g_{m}$ or
$x_m=0$}
\end{cases}
\end{multline*}
where $[g_0: \, \cdots \,  :g_m]$ is the equivalence class of the sequence
$(\seq{g})\in G^{m+1}$ by the equivalence relation
\[ (\seq{g}) \sim (gg_0, gg_1, \, \ldots \,  ,gg_m)     \]
for any $g\in G$.

In the non-homogeneous coordinates on $(BG)_n$ the above relations
take the form
\begin{multline*}
(t_1, \, \ldots \,  , t_m, [h_1| \, \cdots \,  |h_m]) \sim\\
\sim \begin{cases}
(t_2, \, \ldots \,  , t_m, [h_2| \, \cdots \,  |h_m])&\text{for $t_1=0$ or
$h_0=e$} \\
(t_1, \, \ldots \,  , \hat{t_i}, \dots , t_m, [h_1| \, \cdots \,  |h_ih_{i+1}|
\, \cdots \,  |h_m])&
\text{for $t_{i}=t_{i+1}$ or $h_i=e$} \\
(t_1, \, \ldots \,  , t_{m-1}, [h_1| \, \cdots \,  |h_{m-1}])& \text{for
$t_m=1$
or $h_m=e$}
\end{cases}
\end{multline*}

The equivalence class of a sequence $(\seq[m]{x},[g_0: \, \cdots \,  :g_m] )$
will be
denoted by $\vl \seq[m]{x}, [g_0: \, \cdots \,  :g_m] \vr$ and the equivalence
class
of a sequence $(t_1, \, \ldots \,  , t_m, [h_1| \, \cdots \,  |h_m])$ will be
denoted by $\vl t_1, \, \ldots \,  , t_m, [h_1| \, \cdots \,  |h_m]\vr$.

The space $BG$ is the quotient of the disjoint union
$\coprod\limits_{m=0}^\infty \gD^m\times G^{m}$
by the above equivalence relations.

The projection $EG \ra BG$ is given by the formula
\[ \vl \seq[m]{x}, \seq[m]{g} \vr \ras \vl \seq[m]{x}, [g_0: \, \cdots \,
:g_m] \vr\]
or in the non-homogeneous coordinates by
\[\vl t_1, \, \ldots \,  , t_m, h_0[h_1| \, \cdots \,  |h_m]\vr \ras
\vl t_1, \, \ldots \,  , t_m, [h_1| \, \cdots \,  |h_m]\vr      \]

Sometimes it is convenient to write the elements of $EG$ and $BG$ in
the form
\[\vl \seq{m}{x},  h_0[h_1| \, \cdots \,  |h_m]\vr\]
and
\[\vl \seq{m}{x},  [h_1| \, \cdots \,  |h_m]\vr\]
respectively, which is mixture of the  baricentric coordinates on
$\gD^m$ and homogeneous coordinates on $G^{m+1}$ or $G^{m}$.

\subsection{A simpicial description of the geometric bar construction.}
\label{app3}
The above definitions of $EG$ and $BG$ can be interpreted in terms of
geometric realizations of some simplicial objects. Actually, to every
topological group $G$ one can assign simplicial topological groups
$EG.$ and $BG.$ defined as follows.

$EG_n =G^{n+1}$, the face homomorphisms
$\del_i:EG_n\ra EG_{n-1}$ are given  by the formula
\[ \del_i(\seq{g}) = (\seq[i-1]{g} ,\widehat{g_{i}},g_{i+1}, \, \ldots \,
,g_n)\] or in the non-homogeneous coordinates by
\begin{multline*}
\del_i(h_0[h_1| \, \cdots \,  |h_n])=
\begin{cases}
h_0h_1[ h_2|h_3| \, \cdots \,  |h_n] & \text{for} \quad i=0\\{}
h_0[h_1| \, \cdots \,  |h_i\cdot h_{i+1}| \dots |h_n] & \text{for} \quad 0< i
< n\\{}
h_0[h_1| \, \cdots \,  |h_{n-1}]& \text{for} \quad i=n
\end{cases}
\end{multline*}
The degeneracy homomorphism $s_i: EG_n\ra EG_{n+1}$ are
given by the formula
\[ s_i(\seq{g}) = (\seq[i-1]{g}, g_i, g_i, g_{i+1}, \, \ldots \,  ,g_n)\]
or in the non-homogeneous coordinates by
\begin{multline*}
s_i(h_0[h_1| \, \cdots \,  |h_n]) = \begin{cases}
h_0[e|h_1| \, \cdots \,  |h_n] & \text{for} \quad i=0\\{}
h_0[h_1| \, \cdots \,  |h_i| e| h_{i+1}| \dots |h_n] & \text{for} \quad 0< i
< n\\{}
h_0[h_1| \, \cdots \,  |h_n|e]& \text{for} \quad i=n
\end{cases}
\end{multline*}

$EG_n =G^{n}$, the face homomorphisms $\del_i:BG_n\ra
BG_{n-1}$ in the homogeneous coordinates on $G^{n} = G^{n+1}/G$ is
given by the formula
\[\del_i([g_0: \, \cdots \,  : g_n]) = [g_0: \, \cdots \,  g_{i-1}:
\widehat{g_i}:\, \cdots \,
:g_n]\]
or in the non-homogeneous coordinates by
\begin{multline*}
 \del_i([h_1| \, \cdots \,  |h_n]) =
\begin{cases}
[h_2|h_3| \, \cdots \,  |h_n] & \text{for} \quad i=0\\{}
[h_1| \, \cdots \,  |h_i\cdot h_{i+1}| \dots |h_n] & \text{for} \quad 0< i
< n\\{}
[h_1| \, \cdots \,  |h_{n-1}]& \text{for} \quad i=n
\end{cases}
\end{multline*}

The degeneracy homomorphisms $s_i: BG_n\ra BG_{n+1}$ are
given by the formula
\[ s_i([g_0: \, \cdots \,  : g_n]) = [g_0: \, \cdots \,  g_{i-1}: g_i: g_i :\,
\cdots \,
:g_n]\]
or in the non-homogeneous coordinates by
\begin{multline*}
s_i([h_1| \, \cdots \,  |h_n]) =
\begin{cases}
[e|h_1| \, \cdots \,  |h_n] & \text{for} \quad i=0\\{}
[h_1| \, \cdots \,  |h_i| e| h_{i+1}| \dots |h_n] & \text{for} \quad 0< i
< n\\{}
[h_1| \, \cdots \,  |h_n|e]& \text{for} \quad i=n
\end{cases}
\end{multline*}

The geometric realization $|EG.|$ of the simplicial space
$EG.$ is by definition the quotient space of the infinite
disjoint union $\coprod\limits_{n=0}^\infty \gD^n\times G^{n+1}$
by the equivalence relations
\begin{gather*}
(\del^ix,\bar{g}) \sim (x,\del_i\bar{g}) \quad \text{for}\quad
(x,\bar{g})\in \gD^{n-1}\times G^{n+1}\\
(s^ix,\bar{g}) \sim (x,s_i\bar{g}) \quad \text{for}\quad
(x,\bar{g})\in \gD^{n+1}\times G^{n+1}
\end{gather*}
where the maps $\del^i : \gD^{n-1} \ra \gD^n$ and $s^i: \gD^{n+1}
\ra \gD^n$ are defined in the baricentric coordinates by
\begin{align*}
\del^i(\seq{x}) &= (\seq[i-1]{x},0,x_{i},\, \ldots \,  ,x_n)\\
s^i(\seq{x}) &= (\seq[i-1]{x},x_i+x_{i+1},x_{i+2},\, \ldots \,  ,x_n)
\end{align*}
and in the non-homogeneous coordinates by
\begin{align*}
\del^i(\seq[n+1]{t}) &= (\seq[i]{t},t_i,t_{i+1},\, \ldots \,  ,t_{n+1})\\
s^i(\seq[n+1]{t}) &= (\seq[i]{t},\widehat{t_{i+1}},t_{i+2},\, \ldots \,
,t_{n+1})
\end{align*}

Similarly, the geometric realization $|BG.|$ of the simplicial space
$BG.$ is the quotient space of the disjoint union
$\coprod\limits_{n=0}^\infty\gD^n\times G^{n}$
by the equivalence relations
\begin{gather*}
(\del^ix,\bar{g}) \sim (x,\del_i\bar{g}) \quad \text{for}\quad
(x,\bar{g})\in \gD^{n-1}\times G^{n}\\{}
(s^ix,\bar{g}) \sim (x,s_i\bar{g}) \quad \text{for}\quad
(x,\bar{g})\in \gD^{n+1}\times G^{n}
\end{gather*}

\subsection{Group structures on $EG$ and $BG$.}

The usefulness of the non-homogeneous coordinates $t_1, \, \ldots \,  ,
t_{n}$ on $\gD^n$ comes from the fact that they supply a very simple
formula
\begin{gather*}
(t_1, \, \ldots \,  , t_{n}) \times (t_{n+1}, \, \ldots \,  ,t_{n+m+1}) \ras
(t_{\gs(1)}, t_{\gs(2)}, \, \ldots \,  , t_{\gs(n+m+1)})
\end{gather*}
for a homeomorphism pairing
\[\gD^n \times \gD^{m} \lra \gD^{n+m}   \]
where \gs is a permutation of the set $\{1,2, \, \ldots \,  ,n+m+1\}$ such
that
\[ t_{\gs(1)} \leq t_{\gs(2)} \leq \, \ldots \,  \leq t_{\gs(n+m+1)}.\]

Using this pairing we can define, for $G$ an abelian topological
group,  commutative, associative, and continuous pairings
\begin{multline*}
\vl t_1, \, \ldots \,  , t_{n} , h[h_1| \, \cdots \,  |h_n]\vr +
\vl t_{n+1},\, \ldots \,  ,t_{n+m+1}, h'[h_{n+1}|\, \cdots \,  |h_{n+m+1}]\vr
=\\
= \vl t_{\gs(1)}, \, \ldots \,  ,t_{\gs(n+m)}, h\cdot h'[h_{\gs(1)}| \, \cdots
\,  |
h_{\gs(n+m+1)}]\vr
\end{multline*}
and
\begin{multline*}
\vl t_1, \, \ldots \,  , t_{n} , [h_1| \, \cdots \,  |h_n]\vr +
\vl t_{n+1},\, \ldots \,  ,t_{n+m+1}, [h_{n+1}|\, \cdots \,  |h_{n+m+1}]\vr =\\
= \vl t_{\gs(1)}, \, \ldots \,  ,t_{\gs(n+m+1)}, [h_{\gs(1)}| \, \cdots \,  |
h_{\gs(n+m+1)}]\vr
\end{multline*}
\n on $EG$ and $BG$ respectively, which induce group structure on
these spaces \cite{milg-bar}.

The above group parings can be  interpreted as the compositions
\begin{gather*}
EG\times EG = |EG.|\times |EG.|
\overset{\bar{\gf}}{\lra}
|EG.\times EG.| \overset{\bar{\gp}}{\lra}
|E(G\times G).| \overset{|\bar{\mu}|}{\lra}|EG.|=EG\\
BG\times BG = |BG.|\times |BG.|\overset{\gf}{\lra}
|BG.\times BG.| \overset{\gp}{\lra}
|B(G\times G).| \overset{|\mu|}{\lra}|BG.|=BG
\end{gather*}
where $\bar{\gf}$ and $\gf$ are the commutativity of geometric
realization and product operations homeomorphisms, $\bar{\gp}$ and
$\gp$ are induced by the maps
\begin{gather*}
EG_n \times EG_n \lra EG_n\\
(\seq{g})\times (\seq{g'}) \ras ((g_0,g'_0), \, \ldots \,  ,(g_n,g'_n))\\
\intertext{and}
BG_n \times BG_n \lra BG_n\\
[g_0:\, \cdots\, :g_n]\times [g'_0:\, \cdots\, :g'_n]
\ras [(g_0,g'_0):\, \cdots\, :(g_n,g'_n)]
\end{gather*}
respectively, and $|\bar{\mu}|, |\mu|$ are induced by the simplicial
morphisms
\begin{gather*}
\bar{\mu}: E(G\times G). \lra EG.\\
\bar{\mu}((g_0,g'_0), \, \ldots \,  ,(g_n,g'_n)) = (g_0g'_0, \, \ldots \,
,g_ng'_n)\\
\intertext{and}
\mu: B(G\times G). \lra BG.\\
\mu([(g_0,g'_0):\, \cdots\, :(g_n,g'_n)]) = [g_0g'_0:\, \cdots\,
:g_ng'_n]
\end{gather*}
which are well defined only when the multiplication pairing
$\mu:G\times G \ra G$ is a homomorphism, or equivalently, when $G$ is
an abelian group.

\end{document}